\documentclass[11pt]{amsart}
\usepackage{amsmath}
\usepackage{amscd}
\usepackage{pb-diagram}
\usepackage{comment}
\usepackage{graphicx}
\usepackage{amssymb}
\usepackage{lmodern}

\def\refeq#1{\if\workingver y(\ref{#1})-[[#1]]\else(\ref{#1})\fi}
\def\refth#1{\if\workingver y\ref{#1}-[[#1]]\else\ref{#1}\fi}
\def\mylabel#1{\if\workingver y\label{#1}{\bf\ \ [[#1]]\ \ }\else\label{#1}\fi}
\def\mybibitem#1{\if\workingver y\bibitem{#1}{\bf\ \ [[#1]]\ \
}\else\bibitem{#1}\fi}

\newfont{\msam}{msam10}
\newfont{\msbm}{msbm10}
\newcommand{\mvmap}{{\,\overrightarrow{\to}\,}}

\def\articletheorems{
\newtheorem{thm}{Theorem}[section]
\newtheorem{lem}[thm]{Lemma}

\newtheorem{defn}[thm]{Definition}
\newtheorem{cor}[thm]{Corollary}
\newtheorem{prop}[thm]{Proposition}

\newtheorem{algo}{Algorithm}[section] 
\newtheorem{alg}[thm]{Algorithm}  

}

\def\mapright#1{\stackrel{#1}{\rightarrow}}

\def\cB{\text{$\mathcal B$}}

\def\cM{\text{$\mathcal M$}}

\def\cP{\text{$\mathcal P$}}

\def\cT{\text{$\mathcal T$}}

\def\cV{\text{$\mathcal V$}}
\def\cW{\text{$\mathcal W$}}

\newcommand{\id}{\operatorname{id}}
\newcommand{\cl}{\operatorname{cl}}
\newcommand{\opn}{\operatorname{opn}}

\newcommand{\dom}{\operatorname{dom}}

\newcommand{\card}{\operatorname{card}}

\newcommand{\sgn}{\operatorname{sgn}}

\newcommand{\im}{\operatorname{im}}

\renewcommand{\emptyset}{\varnothing}

\newcommand{\Con}{\operatorname{Con}}
\newcommand{\Inv}{\operatorname{Inv}}

\def\proof{{\bf Proof:\ }}

\def\begeq#1{\begin{equation}\mylabel{#1}}
\def\endeq{\end{equation}}

\def\mathobj#1{\mbox{$#1$}}

\def\PP{\mathobj{\mathbb{P}}}

\def\RR{\mathobj{\mathbb{R}}}

\def\ZZ{\mathobj{\mathbb{Z}}}

\def\implies{\;\Rightarrow\;}

\def\setof#1{\mbox{$\{\,#1\,\}$}}

\newcommand{\bdy}{{\bf \partial}}

\def\0#1{\hbox{\kern25pt}$ #1 $\\}
\def\1#1{\hbox{\kern40pt}$ #1 $\\}
\def\2#1{\hbox{\kern55pt}$ #1 $\\}
\def\3#1{\hbox{\kern70pt}$ #1 $\\}

\newcounter{li}

\def\begalg#1{\begin{algo}\mylabel{#1}\normalshape:\small\baselineskip 10pt\\}
\def\endalg{\end{algo}}

\def\Figures(include=#1,cat=#2){
  \renewcommand{\textfraction}{.20}
  \renewcommand{\topfraction}{.80}
  \renewcommand{\bottomfraction}{.80}
  \renewcommand{\floatpagefraction}{.80}
  \newcount\figcount
  \figcount=0
  \let\includefigures=#1
  \def\figcat{#2}
}

\def\FigureFromFile[#1][#2](#3)#4
{
  \begin{figure}[htbp]
     \global\advance\figcount by 1
     \if\includefigures y\special{anisoscale #1.wmf, \the\hsize #2}\fi
     \vspace{#2}
     \caption{#4}
     \mylabel{#3}
   \end{figure}
}

\def\FigureFromFileTwoD[#1][#2,#3](#4)#5
{
  \begin{figure}[htbp]
     \global\advance\figcount by 1
     \if\includefigures y\special{anisoscale #1.wmf, #2 #3}\fi
     \vspace{#2}
     \caption{#5}
     \mylabel{#4}
   \end{figure}
}

\def\FigureF<#1>[#2](#3)#4
{
  \begin{figure}[htbp]
     \global\advance\figcount by 1
     \if\includefigures y\special{anisoscale \figcat/fig\number\figcount.wmf,
       \the\hsize #2}
     \fi
     \if\includefigures p
       \leavevmode
       \epsfxsize=\hsize
       \epsffile{#1}
     \fi
     \if\includefigures y
          \vspace{#2}
     \fi
     \caption{#4}
     \mylabel{#3}
   \end{figure}
}

\def\Figure[#1](#2)#3
{
  \begin{figure}[htbp]
     \global\advance\figcount by 1
     \if\includefigures y\special{anisoscale \figcat/fig\number\figcount.wmf,
       \the\hsize #1}
     \fi
     \if\includefigures p
       \leavevmode
       \epsfxsize=\hsize
       \epsffile{fig\number\figcount.eps}
     \fi
     \if\includefigures y
          \vspace{#1}
     \fi
     \caption{#3}
     \mylabel{#2}
   \end{figure}
}

\def\0{\hbox{\kern5pt}}
\def\1{\hbox{\kern20pt}}
\def\2{\hbox{\kern35pt}}
\def\3{\hbox{\kern50pt}}
\def\4{\hbox{\kern65pt}}
\def\5{\hbox{\kern80pt}}
\def\6{\hbox{\kern95pt}}

\def\kw#1{\textbf{#1}}

\def\kwif{\kw{if}\;}
\def\kwthen{\;\kw{then}\;}
\def\kwelse{\kw{else}\;}
\def\kwendif{\kw{endif}}

\def\kwforeach{\kw{foreach}\;}

\def\kwfunction{\kw{function}\;}

\def\kwreturn{\kw{return}\;}
\def\kwdo{\;\kw{do}\;}
\def\kwenddo{\kw{enddo}}

\def\kwcontinue{\kw{continue}}

\def\kwand{\;\kw{and}\;}

\def\kwnot{\;\kw{not}\;}

\def\Arg{\mbox{$\operatorname{arg}\,$}}

\renewcommand{\mvmap}{{\,\overrightarrow{\to}\,}}

\newcommand{\rank}{\operatorname{rank}}

\newcommand{\mouth}{\operatorname{mo}}

\newcommand{\CR}{\operatorname{CR}}
\newcommand{\Sol}{\operatorname{Sol}}

\def\adhl{\prec}

\def\mapright#1{\stackrel{#1}{\longrightarrow}}

\def\E#1{\hat{#1}}

\def\regclass#1{\langle #1 \rangle_{\cV}}
\def\vclass#1{[#1]_{\cV}}
\def\class#1{[#1]}

\def\leqV{\leq_{\cV}}

\def\rpthV{\rightarrow_{\cV}}
\def\epthV{\leftrightarrow_{\cV}}

\def\rpthVA#1{\stackrel{#1}{\rightarrow}_{\cV}}
\def\epthVA#1{\stackrel{#1}{\leftrightarrow}_{\cV}}

\def\rpthA#1{\stackrel{#1}{\rightarrow}}
\def\epthA#1{\stackrel{#1}{\leftrightarrow}}

\renewcommand{\subset}{\subseteq}

\newcommand{\pto}{{\nrightarrow}}

\articletheorems

\title[Conley-Morse-Forman theory]
{
  Conley-Morse-Forman theory for combinatorial multivector fields
}
\author{Marian Mrozek}
\thanks{}
\address{Marian Mrozek, Division of Computational Mathematics,
  Institute of Computer Science and Computational Mathematics,
  Faculty of Mathematics and Computer Science,
  Jagiellonian University, ul.~St. \L{}ojasiewicza 6, 30-348~Kra\-k\'ow, Poland.
}
\email{Marian.Mrozek@ii.uj.edu.pl}
\thanks{This research is partially supported
       by the Polish National Science Center under Ma\-estro Grant No. 2014/14/A/ST1/00453.}
\subjclass[2010]{Primary
54H20, 
37B30, 
37B35, 
65P99  
; Secondary
57Q10, 
18G35, 
55U15, 
06A06  
}
\keywords{combinatorial vector field, Conley index theory, discrete Morse theory, attractor, repeller,
Morse decomposition, Morse inequalities, topology of finite sets}

\begin{document}

\begin{abstract}
We introduce combinatorial multivector fields, associate with them multivalued dynamics
and study their topological features.
Our combinatorial multivector fields generalize combinatorial vector fields of Forman.
We define isolated invariant sets, Conley index, attractors, repellers and Morse decompositions.
We provide a topological characterization of attractors and repellers and prove Morse inequalities.
The generalization aims at algorithmic analysis of dynamical systems through combinatorialization of flows
given by differential equations and through sampling dynamics in physical and numerical experiments.
We provide a prototype algorithm for such applications.
\end{abstract}

\maketitle

\section{Introduction}
\label{sec:intro}

In the late 90's of the twentieth century Robin Forman \cite{Fo98a} introduced the concept
of a combinatorial vector field
and presented a version of Morse theory for acyclic combinatorial vector fields.
In another paper \cite{Fo98b} he studied combinatorial
vector fields without acyclicity assumption, extended the notion of
the chain recurrent set to this setting and proved Conley type
generalization of Morse inequalities.

Conley theory \cite{Co78} is a generalization of Morse theory
to the setting of non-necessarily gradient or gradient-like
flows on locally compact metric spaces.
In this theory the concepts
of a critical point and its  Morse index are replaced by the more general concept of an isolated
invariant set and its Conley index.
The Conley theory reduces to the Morse theory in the case of a flow
on a smooth manifold defined by a smooth gradient vector field
with non-degenerate critical points.

Recently, T. Kaczynski, M. Mrozek and Th. Wanner \cite{KMW}
defined the concept of an isolated invariant set and the Conley index
in the case of a combinatorial vector field on the collection of simplices of a simplicial
complex and observed that such a combinatorial field has a counterpart
on the polytope of the simplicial complex in the form of a multivalued,
upper semicontinuous, acyclic valued and homotopic to identity map.

\begin{figure}
  \includegraphics[width=0.65\textwidth]{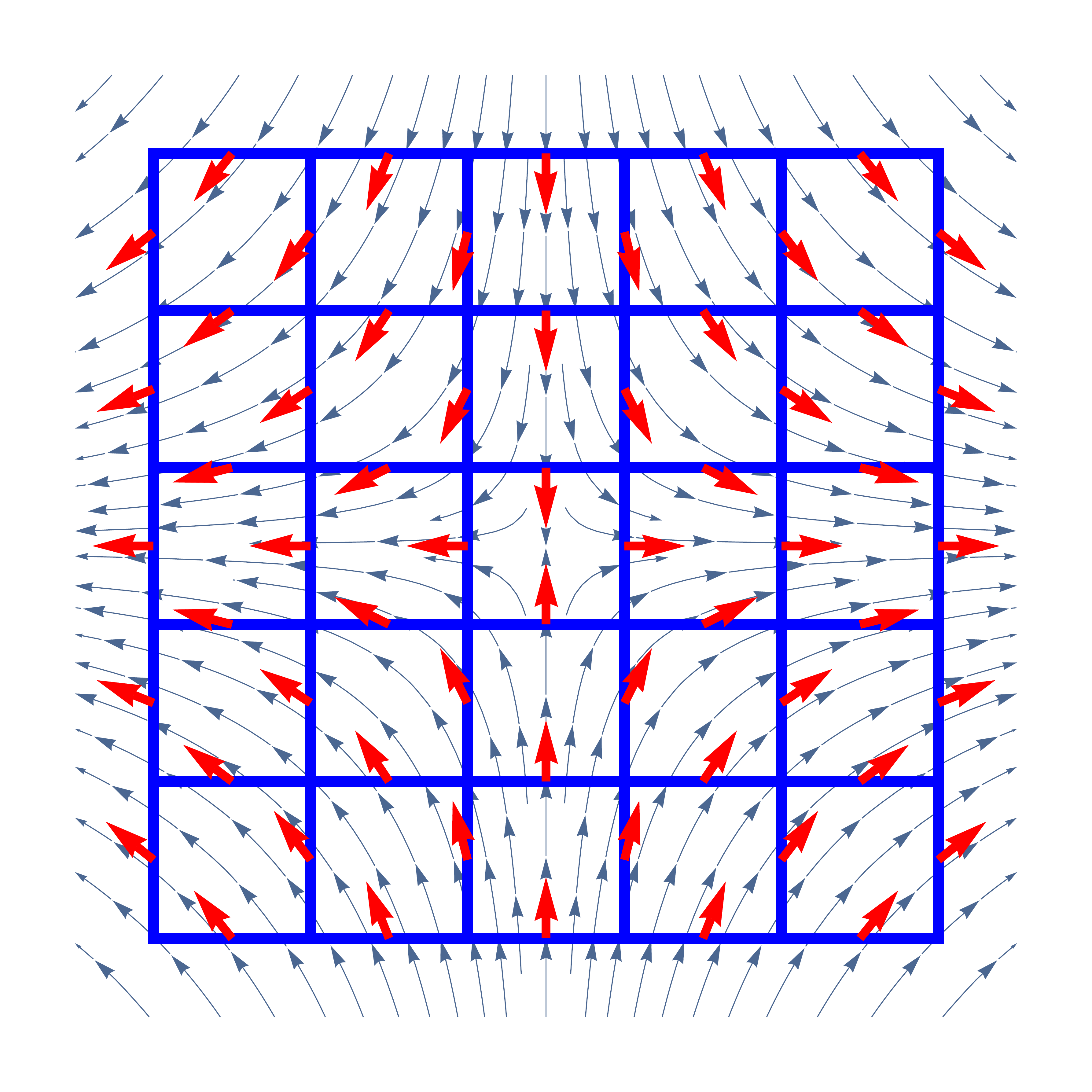}
  \caption{An averaging of a smooth vector field (small arrows)
  along the one-dimensional faces of a cubical grid. Unfortunately, the resulting collection
  of combinatorial vectors (large arrows) does not satisfy the partition requirement of the combinatorial
  vector field of Forman.}
  \label{fig:multivectorsDraws}
\end{figure}

The aim of this paper is to combine the ideas of Forman with some classical
concepts of topological dynamics in order to obtain an algorithmic tool for studying sampled dynamics,
that is dynamics known only via a finite set of samples obtained from a physical or numerical experiment.
The method to achieve this aim is the combinatorialization of classical dynamics.
By this we mean constructing an analogue of classical topological dynamics set up
in finite combinatorial spaces: simplicial complexes, cubical sets (also called cubical complexes) or more generally
cellular complexes. Such spaces are equipped with a natural but non-Hausdorff topology
via the Alexandroff Theorem~\cite{Al1937} and the partial order given by the face relation.
Simplicial complexes in the form of triangular meshes are typically used in visualization of vector fields sampled from data
and the use of topological methods in this field increases \cite{ChMiLa2008,RSHWCS2012,Szy2012}.
In gene regulatory networks a frequent method used to analyse the associated dynamics
is Thomas' formalism \cite{SnTh1993} leading to the study of dynamics on cubical grids \cite{BCRG2004}.
The proposed combinatorialization may also serve as a very concise description of
the qualitative features of classical dynamics.

Forman's vector fields seem to be a natural tool for a concise approximation and description
of the dynamics of differential equations and more generally flows.
For instance, given a cubical grid in $\RR^d$  and a
vector field, it is natural to set up arrows in the combinatorial setting
of the grid by taking averages of the vectors in the vector field along
the codimension one  faces of the grid.
Unfortunately, in most cases
such a procedure does not lead to a well defined combinatorial vector field in the sense
of Forman. This is because in the Forman theory
the combinatorial vectors
together with the critical cells have to constitute a partition.
In particular, each non-critical, top-dimensional cell has to be paired with precisely one
cell in its boundary. Such a requirement is not satisfied by a typical space discretization
of a vector field (see Figure~\ref{fig:multivectorsDraws}).
In order to overcome these limitations introduce and study combinatorial {\em multivector} fields,
a generalization of Forman's combinatorial vector fields.
Similar but different generalizations of Forman's combinatorial vector fields are proposed by
Wisniewski and Larsen \cite{Wisn2008} in the study of piecewise affine control systems
and by Freij \cite{Freij2009} in the combinatorial study of equivariant discrete Morse theory.

We extend the concepts of isolated invariant set and Conley index introduced
in \cite{KMW} to combinatorial multivector fields.
We also define attractors, repellers, attractor-repeller pairs and Morse decompositions
and provide a topological characterization of attractors and repellers.
These ideas are novel not only for combinatorial multivector fields
but also for combinatorial vector fields.
Furthermore, we prove the Morse equation  for Morse decompositions.
We deduce from it  Morse inequalities.
They generalize the Morse inequalities proved by Forman in \cite{Fo98b}
for the Morse decomposition consisting of basic sets of a combinatorial vector field
to the case of general Morse decompositions for combinatorial multivector fields.

The construction of the chain complex, an algebraic structure needed in our study, is complicated in the case
of a general cellular complex. This is in contrast to the case of a simplicial complex or a cubical set.
To keep things simple but general, in this paper we work in the algebraic setting of chain complexes
with a distinguished basis, an abstraction of the chain complex of a simplicial, cubical or cellular complex
already studied by S.\ Lefschetz \cite{Le1942}. It is elementary
to see simplicial and cubical complexes  as examples of Lefschetz complex.
A version of Forman theory for combinatorial vector fields on chain complexes
with a distinguished basis was recently proposed by a few authors \cite{JW,Ko,Sk}.
Related work concerns Forman theory on finite topological spaces \cite{Mi2012}.

The organization of the paper is as follows.
In Section~\ref{sec:main} we provide an informal overview of the main results of the paper.
In Section~\ref{sec:examples} we illustrate the new concepts and results with several examples.
In Section~\ref{sec:preli} we gather preliminary definitions and results.
In Section~\ref{sec:multi} we introduce Lefschetz complexes, define combinatorial multivector fields and prove their basic features.
In Section~\ref{sec:solu} we define solutions and invariant sets of combinatorial multivector fields.
In Section~\ref{sec:iso} we study isolated invariant sets of combinatorial multivector fields and their Conley index.
In Section~\ref{sec:attractors-repellers} we investigate attractors,
repellers and attractor-repeller pairs.
In Section~\ref{sec:morse} we define Morse decompositions and prove Morse equation and Morse inequalities.
In Section~\ref{sec:algo} we discuss an algorithm constructing combinatorial multivector fields from clouds of vectors
on the planar integer lattice.
In Section~\ref{sec:extensions} we show a few possible extensions of the theory presented in this paper.
In Section~\ref{sec:conclu} we present conclusions and directions of future research.

\section{Main results.}
\label{sec:main}

In this section we informally present the main ideas and results of the paper.
Precise definitions, statements and proofs will be given in the sequel.

\subsection{Lefschetz complexes.}
\begin{figure}
\begin{center}
  \includegraphics[width=0.65\textwidth]{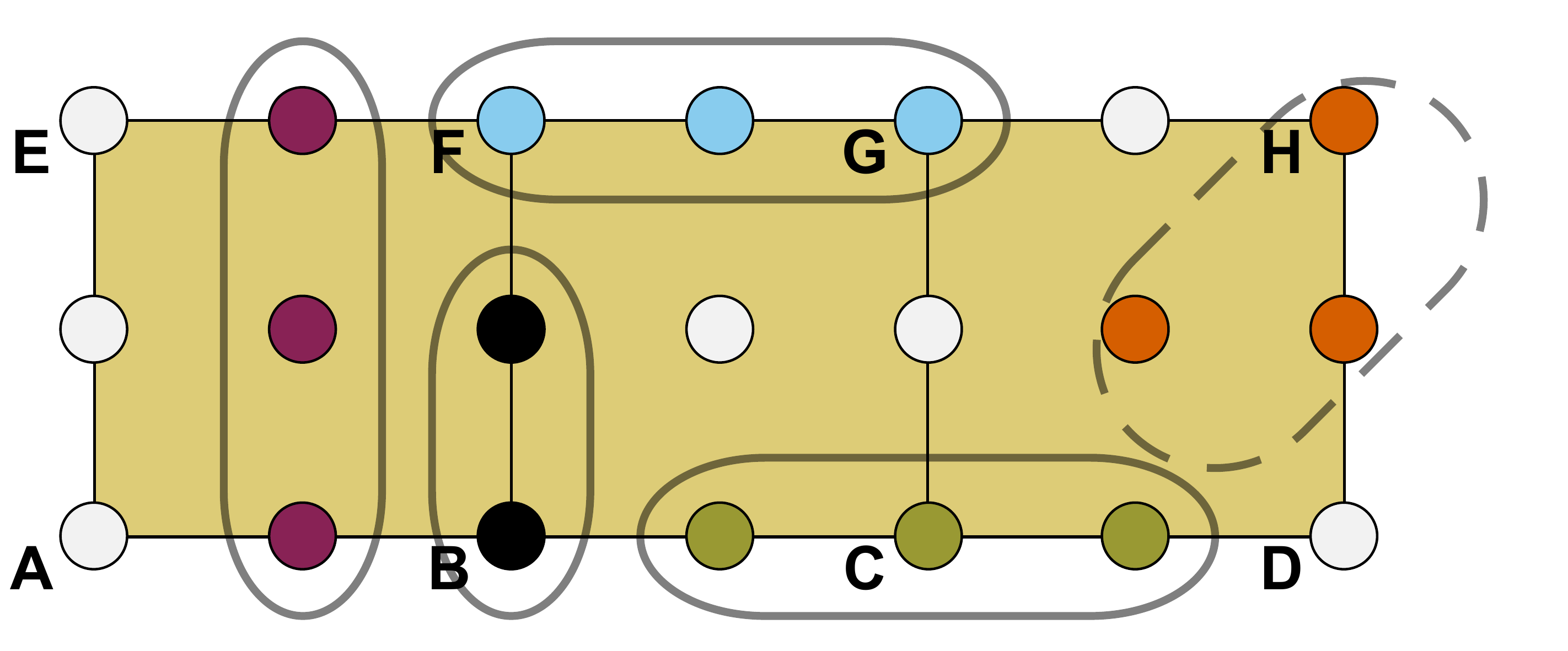}
\end{center}
  \caption{A Lefschetz complex consisting of the collection $K$ of cells of a cubical complex
  with eight vertices ($0$-cells or cells of dimension zero), ten edges ($1$-cells)
  and three squares ($2$-cells).
  Individual cells are marked by a small circle in the center of mass of each cell.
  Four Lefschetz complexes obtained as proper subsets of $K$ are indicated by solid ovals.
  The collection of three cells marked by a dashed oval is not a Lefschetz complex, because
  it is not proper.}
  \label{fig:cw-complex}
\end{figure}

Informally speaking, a {\em Lefschetz complex}
is an algebraization of the simplicial, cubical or cellular complex.
It consists of a finite collection of cells $X$ graded by dimension
and the {\em incidence coefficient} $\kappa(x,y)$ encoding the incidence relation
between cells $x,y\in X$ (see Section~\ref{sec:lefschetz} for a precise definition).
A non-zero value of $\kappa(x,y)$
indicates that the cell $y$ is in the boundary of the cell $x$ and the dimension
of $y$ is one less than the dimension of $x$. The cell $y$ is then called a {\em facet} of $x$.
The family $K$ of all simplices of a simplicial complex \cite[Definition 11.8]{KMW2004},
all elementary cubes in a cubical set \cite[Definition 2.9]{KMW2004}
or, more generally, cells of a cellular complex (finite CW complex, see \cite[Section IX.3]{Ma1991})
are examples of Lefschetz complexes.
In this case the incidence coefficient is obtained from the
boundary homomorphism of the associated simplicial, cubical or cellular chain complex.
A sample Lefschetz complex is presented in Figure~\ref{fig:cw-complex}.
It consists of eight vertices ($0$-cells or cells of dimension zero), ten edges ($1$-cells)
and three squares ($2$-cells).

Condition \eqref{eq:kappa-condition} presented in Section~\ref{sec:lefschetz}
guarantees that
the free group spanned by $X$ together with the linear map given by
$
    \bdy x:=\sum_{y\in X}\kappa(x,y)y
$
is a free chain complex with $X$ as a basis.
By the {\em Lefschetz homology} of $X$ we mean the homology of this chain complex.
We denote it by $H^\kappa(X)$.
In the case of a Lefschetz complex given
as a cellular complex condition \eqref{eq:kappa-condition} is satisfied and the resulting chain complex
and homology is precisely the cellular chain complex and the cellular homology.
Given a Lefschetz complex $X$ we write $p_X$ for the respective {\em Poincar\'e polynomial}, that is,
\[
    p_X(t):=\sum_{i=0}^\infty\rank H^\kappa_i(X)t^i.
\]

The {\em closure} of $A\subset X$, denoted $\cl A$, is obtained by
recursively adding to $A$ the facets of
cells in $A$, the facets of the facets of cells in $A$ and so on.
The set $A$ is {\em closed} if $\cl A=A$ and it is {\em open} if $X\setminus A$ is closed.
The terminology is justified, because the open sets indeed form a $T_0$ topology on $X$.
We say that $A$ is {\em proper} if $\mouth A:=\cl A\setminus A$, which we call the {\em mouth} of $A$,
is closed. Proper sets are important for us, because every proper subset of a Lefschetz complex
with incidence coefficients restricted to this subset is also a Lefschetz complex.
Four Lefschetz complexes being proper subsets of a bigger Lefschetz complex are indicated
in Figure~\ref{fig:cw-complex} by solid ovals.
In the case of a proper subset $X$ of a cellular complex $K$ the Lefschetz homology of $X$
is isomorphic to the relative cellular homology
$
H( \cl X, \mouth X).
$

\subsection{Multivector fields.}
A {\em multivector} is a proper $V\subset K$ such that $V\subset\cl V^\star$ for a unique $V^\star\in V$.
Out of the four examples of Lefschetz complexes in Figure~\ref{fig:cw-complex} only the one enclosing vertex $C$
is not a multivector. A multivector is {\em critical} if $H^\kappa(V)\neq 0$. Otherwise it is {\em regular}.
The only example of a regular multivector in Figure~\ref{fig:cw-complex} is the one enclosing vertex $B$.
Roughly speaking, a regular multivector indicates that $\cl V$, the closure of $V$, may be collapsed to $\mouth V$,
the mouth of $V$. In the dynamical sense this means that  everything flows through $V$. A critical multivector indicates the contrary:
$\cl V$ may not be collapsed to $\mouth V$ and, in the dynamical sense, something must stay inside $V$.

\begin{figure}
\begin{center}
  \hspace*{6pt}\includegraphics[width=0.62\textwidth]{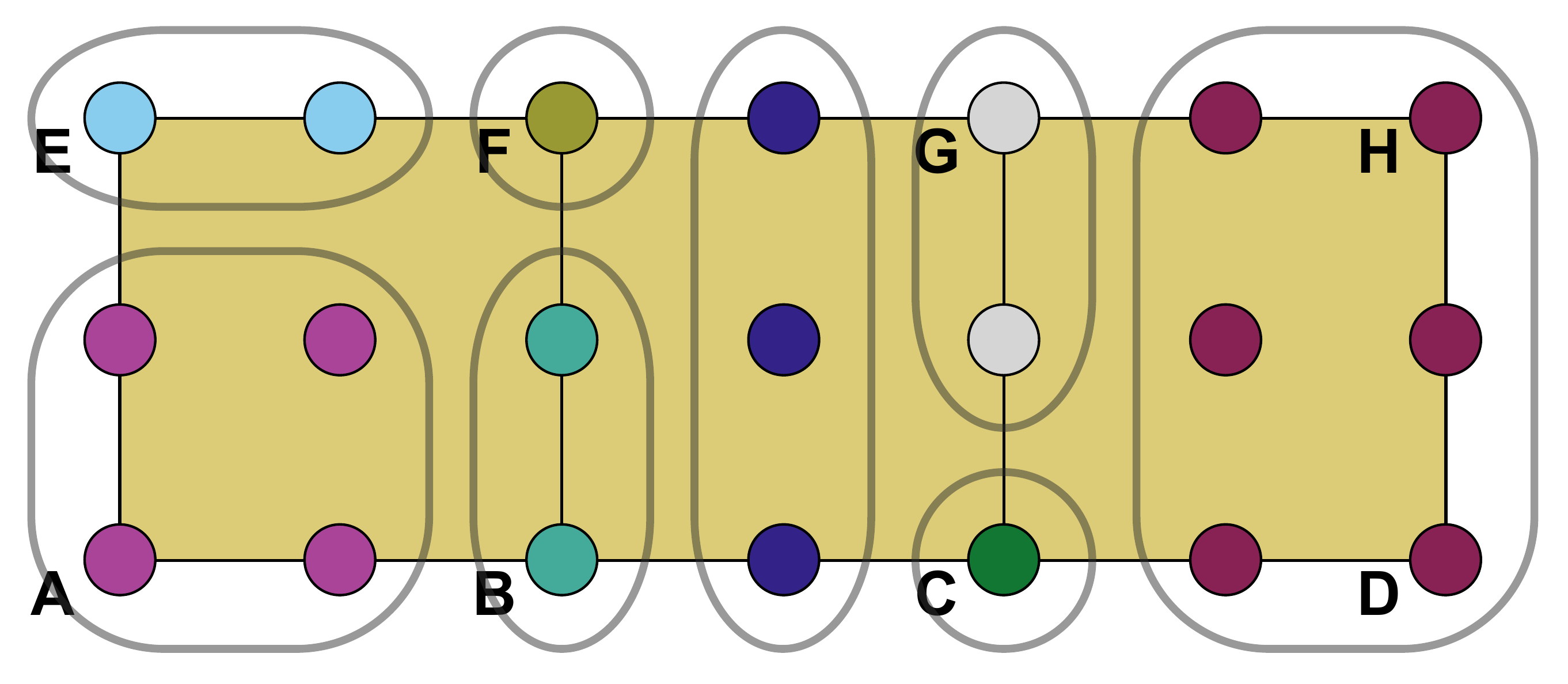}\\
  \includegraphics[width=0.6\textwidth]{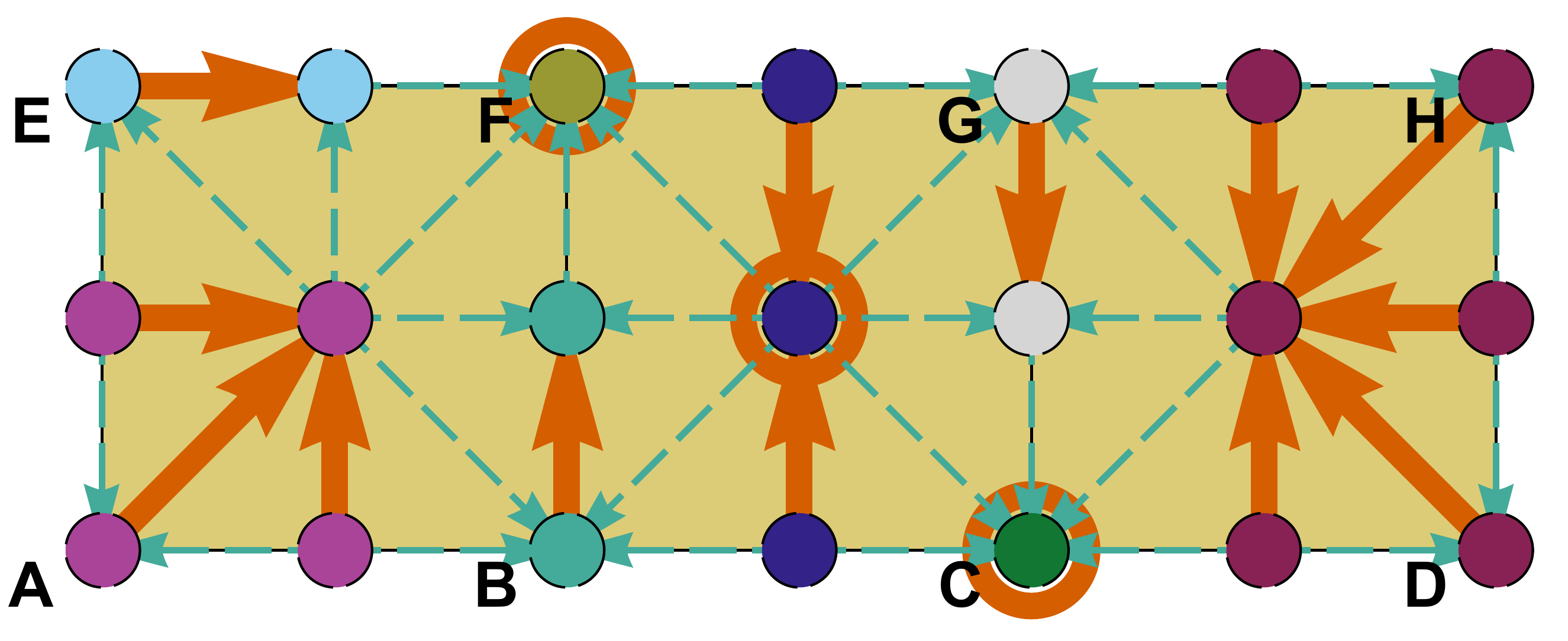}
\end{center}
  \caption{A generalized multivector field as a partition of a Lefschetz complex (top)
  and the associated  directed graph $G_\cV$ (bottom).
  The up-arrows and loops are marked by thick solid lines.
  The down-arrows are marked by thin dashed lines.
  The critical multivectors are $\{F\}$, $\{C\}$ and $\{BC,FG,BCGF\}$.
  }
  \label{fig:MVFpartition}
\end{figure}

A {\em combinatorial multivector field} is a partition $\cV$ of $X$ into multivectors.
We associate dynamics with $\cV$ via a directed graph $G_\cV$ with vertices in $X$
and three types of arrows: {\em up-arrows}, {\em down-arrows} and {\em loops}. Up-arrows have
heads in $V^\star$ and tails in  all the other cells of $V$.
Down-arrows have tails in  $V^\star$ and heads in $\mouth V$.
Loops join $V^\star$ with itself for all critical multivectors $V$.
A sample multivector field is presented in Figure~\ref{fig:MVFpartition}(top)
together with the associated directed graph $G_\cV$ (bottom).
The terminology 'up-arrows' and 'down-arrows' comes from the fact that the dimensions of cells
are increasing along up-arrows and decreasing along down-arrows.
Notice that the up-arrows sharing the same head uniquely determine a multivector.
Therefore, it is convenient to draw a multivector field not as a partition but by marking all up-arrows.
For convenience, we also mark the loops, but the down-arrows are implicit and are usually omitted to keep the drawings
simple.

A multivector may consist of one, two or more cells.
If there are no more than two cells, we say that the multivector is a {\em vector}.
Otherwise we call it a {\em strict} multivector.
Note that the combinatorial multivector field in Figure~\ref{fig:MVFpartition}
has three strict multivectors:
\[
 \{ABFE,AB,AE,A\}, \{BCGF,BC,FG\},   \{CDHG,CD,DH,GH,D,H\}.
\]
Observe that a combinatorial multivector field
with no strict multivectors corresponds to the combinatorial vector field in the sense of Forman \cite{Fo98b}.

A cell $x\in X$ is {\em critical} with respect to $\cV$ if $x=V^\star$ for a critical multivector $V\in\cV$.
A critical cell $x$ is {\em non-degenerate} if the Lefschetz homology of its multivector is zero in all dimensions except one in which it
is isomorphic to the ring of coefficients. This dimension is then the {\em Morse index} of the critical cell.
The combinatorial multivector field in Figure~\ref{fig:MVFpartition} has three critical cells: $F$, $C$ and $BCGF$.
They are all nondegenerate. The cells $F$ and $C$ have Morse index equal zero. The cell $BCGF$ has Morse index equal one.

A {\em solution} of $\cV$ (also called a {\em trajectory} or a {\em walk}) is a bi-infinite,
backward infinite, forward infinite, or finite sequence
of cells such that any two consecutive cells in the sequence form an arrow in the graph  $G_\cV$.
The solution is {\em full} if it is bi-infinite. A finite solution is also called a {\em path}.
The full solution is {\em periodic} if the sequence is periodic. It is {\em stationary} if the sequence is constant.
By the {\em dynamics} of $\cV$ we mean the collection of all solutions.
The dynamics is multivalued in the sense that there may be many different solutions
going through a given cell.

\subsection{Isolated invariant sets.}
Let $X$ be a Lefschetz complex and let $\cV$ be a combinatorial multivector field on $X$.
Assume $S\subset X$ is $\cV$-{\em compatible}, that is, $S$ equals the union of multivectors contained in it.
We say that $S$ is {\em invariant} if for every multivector $V\subset S$ there is
a full solution through $V^\star$ in $S$.
The {\em invariant part} of a subset $A\subset X$
is the maximal, $\cV$-compatible invariant subset of $A$.
A path in $\cl S$ is an {\em internal tangency} to $S$ if the  values at the endpoints of the path
are in $S$ but one of the values is not in $S$.
The set $S$ is {\em isolated invariant} if it is invariant and admits no internal tangencies.

An isolated invariant set is an {\em attractor}, respectively a {\em repeller},
if there is no full solution crossing it which goes away from it in forward, respectively backward, time.
The attractors and repellers have the following topological characterization in terms of the $T_0$ topology
of $X$ (see Sec.~\ref{sec:attractors-repellers}, Theorems~\ref{thm:attractor}~and~\ref{thm:repeller}).
\begin{thm}
\label{thm:attractor-repeller-characterization}
An isolated invariant set $S\subset X$ is an attractor, respectively a repeller, if and only if
it is closed, respectively open, in $X$.
\end{thm}

\subsection{Morse inequalities.}
Given a family $\{M_r\}_{r\in{\scriptsize\PP}}$ of mutually disjoint, non-empty, isolated invariant sets,
we define a relation $r\leq r'$ in $\PP$ if there exists a full solution
such that all its sufficiently far terms  belong to $M_r$ and all sufficiently early terms
belong to $M_{r'}$.
Such a full solution is called a  {\em connection} running from $M_{r'}$ to $M_r$.
The connection is {\em heteroclinic} if $r\neq r'$. Otherwise it is called {\em homoclinic}.
We say that $\{M_r\}_{r\in{\scriptsize\PP}}$ is a {\em Morse decomposition} of $X$
if the relation $\leq$ induces a partial order in $\PP$.
The Hasse diagram of this partial order with vertices labelled by the Poincar\'e polynomials
$p_{M_r}(t)$ is called the {\em Conley-Morse} graph of the Morse decomposition
(comp. \cite[Def. 2.11]{BuAtAl2012}).

The Poincar\'e polynomials $p_{M_r}(t)$ are related to the Poincar\'e polynomial $p_X(t)$ via the following theorem
(see  Section~\ref{sec:morse-equation} Theorems~\ref{thm:Morse-equation} and~\ref{thm:Morse-inequalities}).

\begin{thm} (Morse equation and Morse inequalities)
\label{thm:main}
  If $\{M_r\}_{r\in{\scriptsize \PP}}$ is a Morse decomposition of $X$, then
\begin{equation*}
  \sum_{r\in \scriptsize{\PP}}p_{M_r}(t)=p_X(t)+(1+t)q(t)
\end{equation*}
for some polynomial $q(t)$ with nonnegative coefficients.
In particular, for any natural number $k$ we have
\[
   \sum_{r\in{\scriptsize \PP}}\rank H^\kappa_k(M_r)\geq \rank H^\kappa_k(X).
\]
\end{thm}

\begin{figure}
\begin{center}
  \includegraphics[width=0.69\textwidth]{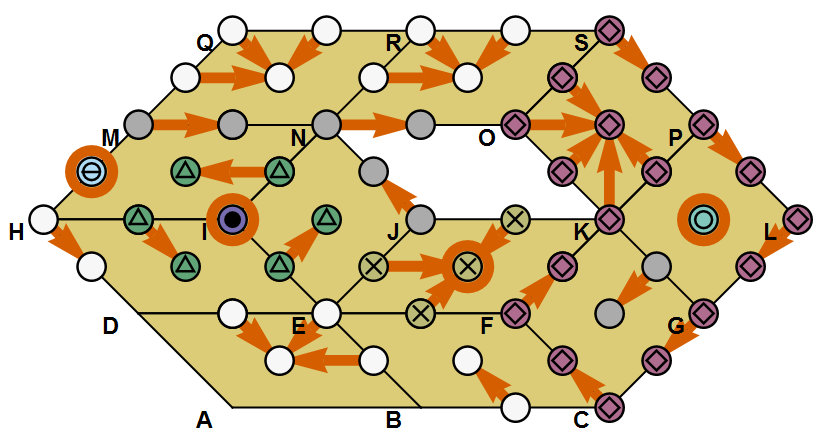}~~~~
  \raisebox{12mm}{\includegraphics[width=0.29\textwidth]{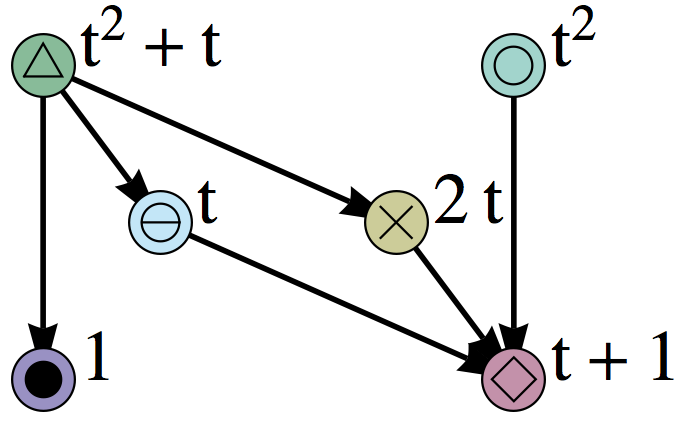}}
\end{center}
  \caption{ A Morse decomposition of a combinatorial multivector field
 (left) and its Conley-Morse graph (right).
  The decomposition consists of six isolated invariant sets.
  Cells in the same sets share the same mark.}
  \label{fig:examplePap1}
\end{figure}

\section{Examples.}
\label{sec:examples}

In this section we present a few examples of combinatorial multivector fields and some of its Morse decompositions.
We begin with
the example in
Figure~\ref{fig:examplePap1}.
Then, we present an example illustrating the
differences between the combinatorial multivector fields and combinatorial vector fields. We complete this section with
examples of combinatorial multivector fields constructed by algorithm CMVF presented in Section~\ref{sec:algo}.
Two of these examples are derived from a planar smooth vector field and one is derived from a cloud of random vectors
on an integer lattice.

\subsection{Attractors and repellers.}
Consider the planar regular cellular complex in Figure~\ref{fig:examplePap1}(left).
It consists of $11$ quadrilaterals and its faces.
A proper subcollection of its $55$ faces, marked by a circle in the center of mass,
forms a Lefschetz complex $X$.
It consists of all cells of the cellular complex except vertices $A$, $B$, $D$ and edges $AB$, $AD$.
A combinatorial  multivector field $\cV$ on $X$ is marked by up-arrows and loops.
The invariant part of $X$ with respect to $\cV$ consists of all cells of $X$ but the cells marked in white.
The Lefschetz homology $H^\kappa(X)\cong H(K,A)$ where $A$ is the cellular complex consisting of vertices $A,B,D$ and edges
$AB$, $AD$. Thus, this is the homology of a pointed annulus. Therefore, $p_X(t)=t$.

Consider the family of six isolated invariant sets
\[
\cM=\setof{M_{\displaystyle{\bullet}},M_{\ominus},M_{\displaystyle{\circ}},M_{\times},M_{\triangle},M_{\Diamond}},
\]
marked in Figure~\ref{fig:examplePap1} with the respective symbols.
The family $\cM$ is a Morse decomposition of $X$.
The respective Poincar\'e polynomials are:
$p_{\displaystyle{\bullet}}(t)=1$, $p_{\ominus}(t)=t$, $p_{\displaystyle{\circ}}(t)=t^2$,
$p_{\times}(t)=2t$, $p_{\triangle}(t)=t^2+t$, $p_{\Diamond}(t)=t+1$.

There are two attractors: stationary $M_{\displaystyle{\bullet}}$ and periodic $M_{\Diamond}$.
There are also two repellers: stationary $M_{\displaystyle{\circ}}$ and periodic $M_{\triangle}$.
The other two isolated invariant sets are neither attractors nor repellers.
The Morse equation takes the form
\[
   2t^2+5t+2=t+(1+t)(2+2t).
\]

\begin{figure}
\begin{center}
  \includegraphics[width=0.45\textwidth]{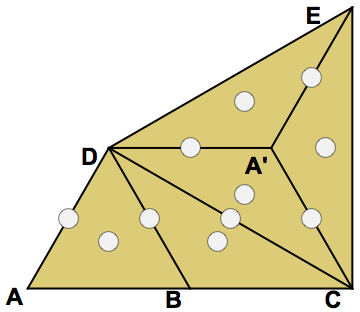}
\end{center}
  \caption{A cellular complex with edges $AD$ and $A'D$ identified. The collection of five triangles and five edges
  marked with a circle in the center of mass is proper, thus forms a Lefschetz complex $X$.
  A direct computation shows that $H^\kappa(X)=0$, hence $X$ is a zero space. However,
  it is easy to see that any combinatorial vector field $\cV$ on $X$
  either has a critical cell or a periodic solution, thus the invariant part
  of $\cV$ is never empty.}
  \label{fig:zero-space-no-acyclic-forman}
\end{figure}

\subsection{Refinements of multivector fields.}
A multivector field $\cW$ is a {\em refinement} of $\cV$ if each multivector in $\cV$
is $\cW$-compatible. The refinement is {\em proper} if the invariant part with respect to $\cW$  of each
regular multivector in $\cV$ is empty.
A {\em Forman refinement} of a multivector field $\cV$
is a vector field $\cW$ which is a proper refinement of $\cV$ such that
each multivector of $\cV$ contains at most one critical vector of $\cW$.
Then $\cW$ has precisely one critical cell
in any critical  multivector of $\cV$.

\begin{figure}
\begin{center}
  \includegraphics[width=0.20\textwidth]{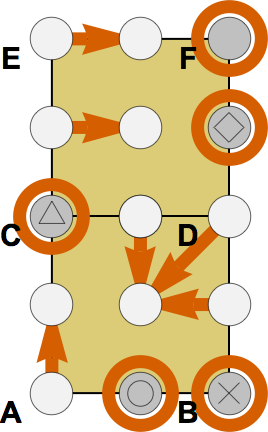}\quad\quad\quad
  \includegraphics[width=0.20\textwidth]{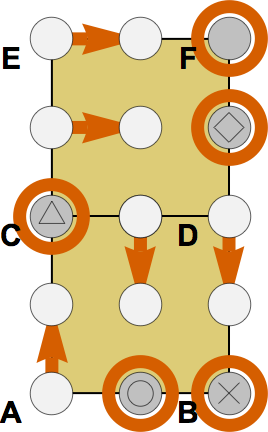}\quad\quad\quad
  \includegraphics[width=0.20\textwidth]{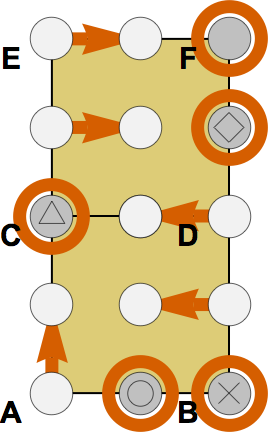}
\end{center}
  \caption{A multivector field (left) and its two different Forman refinements (middle and right).}
  \label{fig:mv-refin}
\end{figure}

The concept of Forman refinement raises two natural questions. The first question is whether a combinatorial multivector
field always admits a Forman refinement. The second question is whether the study of the dynamics
of a combinatorial multivector field which admits a Forman refinement may be reduced to the study
of the dynamics of the refinement.

We do not know what the answer to the first question is.
If the answer is negative, then there exists a multivector field $\cV$ on a Lefschetz complex $X$
such that at least one multivector of $\cV$ cannot be partitioned into vectors with at most one
critical vector in the partition.
There are examples of zero spaces
(Lefschetz complexes with zero homology) which do not admit a combinatorial vector field
with empty invariant part.
They may be constructed by adapting examples of contractible but
not collapsible cellular complexes such as Bing's house \cite{Bi1964} or dunce hat \cite{Zee1964}.
One such example is presented in Figure~\ref{fig:zero-space-no-acyclic-forman}.
This example fulfills all requirements of a multivector except the requirement that
a multivector has precisely one top-dimensional cell, because it has five top-dimensional cells.

Regardless of what is the answer to the first question, even if a given combinatorial multivector field does have a Forman refinement,
in general it is not unique. Figure~\ref{fig:mv-refin} shows a combinatorial multivector field $\cV$ (left)
and its two  different Forman refinements: $\cV_1$ (middle) and $\cV_2$ (right).
The critical cells of all three combinatorial multivector fields are the same:  $AB$, $B$, $C$, $DF$ and $F$.
However, in the case of $\cV$ there are heteroclinic connections running from the critical cell $DF$
to the critical cells  $AB$, $B$ and $C$. In the case of $\cV_1$ there is a heteroclinic connection
running from  $DF$ to $B$ but not to $AB$ nor $C$. In the case of $\cV_2$ there is a heteroclinic connection
running from  $DF$ to $C$ but not to $AB$ nor $B$.
Our next example shows that the differences may be even deeper.
Thus, the answer to the second question is clearly negative.

\begin{figure}
\begin{center}
  \includegraphics[width=0.65\textwidth]{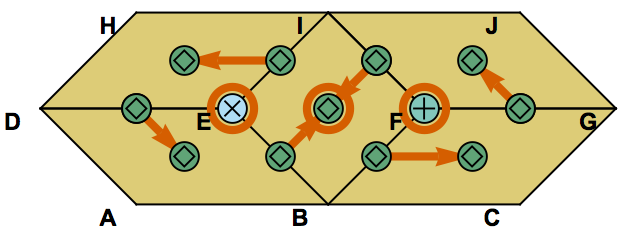}\\
  \includegraphics[width=0.65\textwidth]{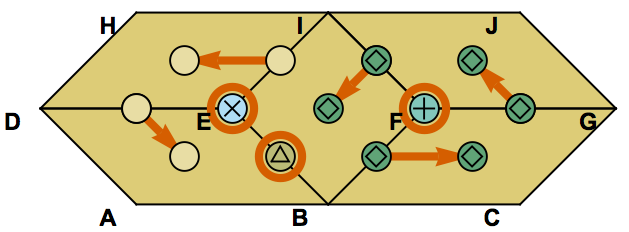}\\
\end{center}
  \caption{A combinatorial multivector field with homoclinic connections and chaotic dynamics (top)
  and one of its two Forman refinements (bottom). Both refinements are deprived of such features.}
  \label{fig:homoclinic}
\end{figure}

\subsection{Homoclinic connections and chaotic dynamics.}
Figure~\ref{fig:homoclinic} presents a combinatorial multivector field $\cV$ on a Lefschetz complex $X$ (top)
and one of its two Forman refinements $\cV_1$ (bottom).
The combinatorial multivector field $\cV$ has
homoclinic connections to the cell $BEIF$. Moreover, it admits chaotic dynamics
in the sense that for each bi-infinite sequence of two symbols marking the two edges $BF$ and $EI$,
there is a full trajectory whose sequence of passing through the edges $BF$ and $EI$ is precisely
the given one. The two Forman refinements of $\cV$ have neither homoclinic connections nor chaotic dynamics.

\begin{figure}
\begin{center}
  \includegraphics[width=0.95\textwidth]{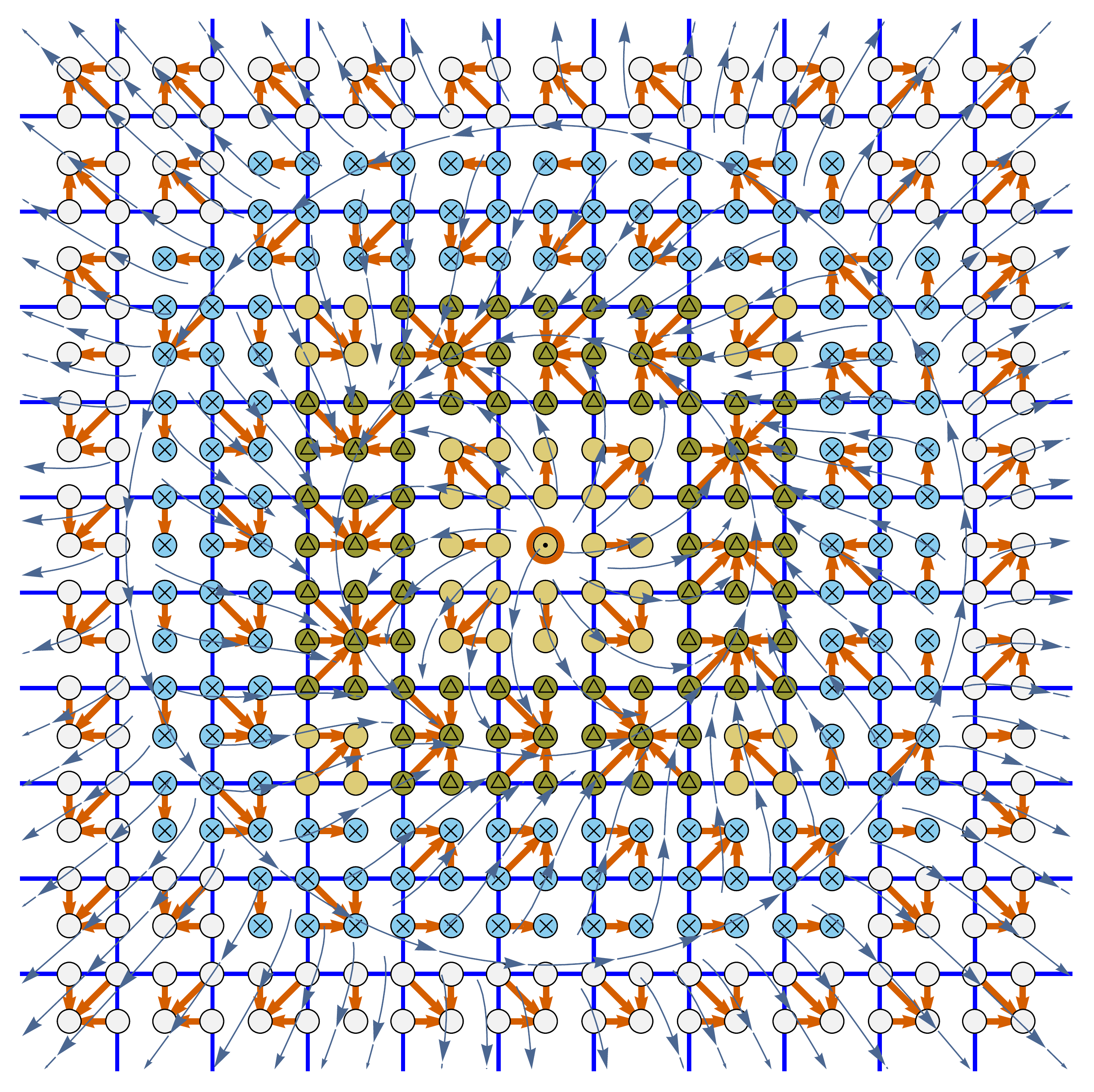}\\
\end{center}
  \caption{A combinatorial multivector field
    modelling the dynamics of the differential equation \eqref{eq:twoCirc-ode}.
    The critical cell in the middle of the grid, marked with a dot, captures the repelling stationary point of \eqref{eq:twoCirc-ode}.
    The isolated invariant set marked with triangles captures the attracting periodic trajectory of  \eqref{eq:twoCirc-ode}.
    The isolated invariant set marked with crosses captures the repelling periodic trajectory of  \eqref{eq:twoCirc-ode}.
    The Conley-Morse graphs of  \eqref{eq:twoCirc-ode} and the combinatorial model coincide.
  }
  \label{fig:AttrRepCircles-cmvf}
\end{figure}

\begin{figure}
\begin{center}
  \includegraphics[width=0.95\textwidth]{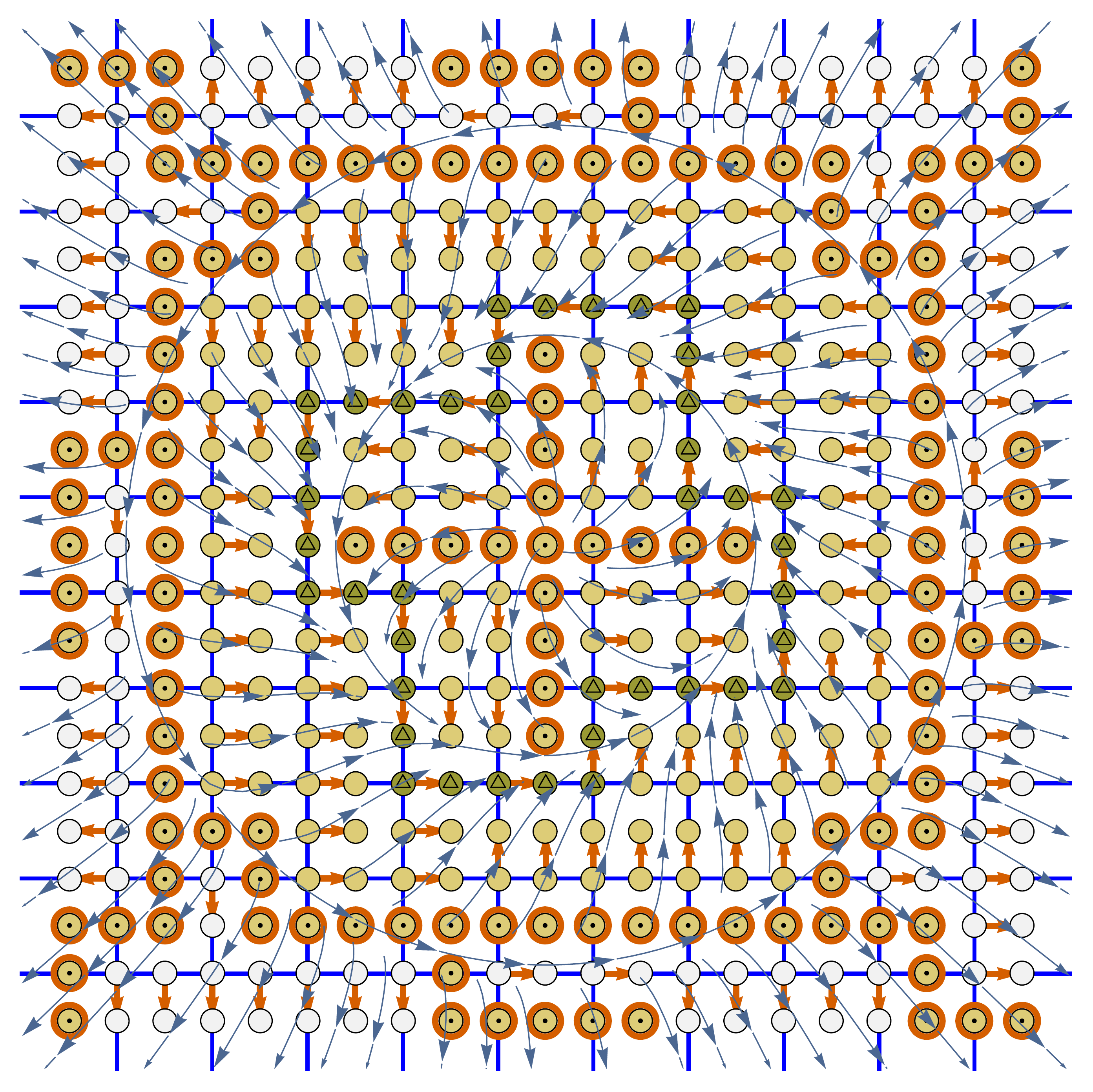}
\end{center}
  \caption{A combinatorial vector field
  modelling the dynamics of the differential equation \eqref{eq:twoCirc-ode}.
  Only the attracting cycle is captured.}
  \label{fig:AttrRepCircles-cvf}
\end{figure}

\subsection{A combinatorial multivector field constructed from a smooth vector field.}
In section~\ref{sec:algo} we present algorithm CMVF. Its input consists of
a collection of classical vectors on an integer, planar lattice.
These may be vectors of a smooth planar vector field evaluated at the lattice points.
However, the algorithm accepts any collection of vectors, also vectors chosen randomly.
It constructs a combinatorial multivector field based on the directions of the classical vectors,
with varying number of strict multivectors: from many to none,
depending on a control parameter.

As our first example consider the vector field of the differential equation
\begin{equation}
\label{eq:twoCirc-ode}
\begin{array}{ccc}
\dot{x}_1 &=&  -x_2 + x_1 (x_1^2 + x_2^2-4) (x_1^2+x_2^2-1) \\
\dot{x}_2 &=&  x_1 + x_2 (x_1^2 + x_2^2-4) (x_1^2+x_2^2-1)
\end{array}
\end{equation}
restricted to the $10\times 10$ lattice of points in the square $[-3,3]\times[-3,3]$.
The equation has three minimal invariant sets: a repelling stationary point at the origin and two invariant circles:
an attracting periodic orbit of radius $1$ and a repelling periodic orbit of radius $2$.
The outcome of algorithm CMVF maximizing the number of strict multivectors is presented in Figure~\ref{fig:AttrRepCircles-cmvf}.
It captures all three minimal invariant sets of \eqref{eq:twoCirc-ode}.
The variant forbidding  strict multivectors  is presented in Figure ~\ref{fig:AttrRepCircles-cvf}.
It captures only the attracting
periodic orbit, whereas the repelling fixed point and repelling periodic trajectory
degenerate into a collection of critical cells.

\begin{figure}
\begin{center}
  \includegraphics[width=0.95\textwidth]{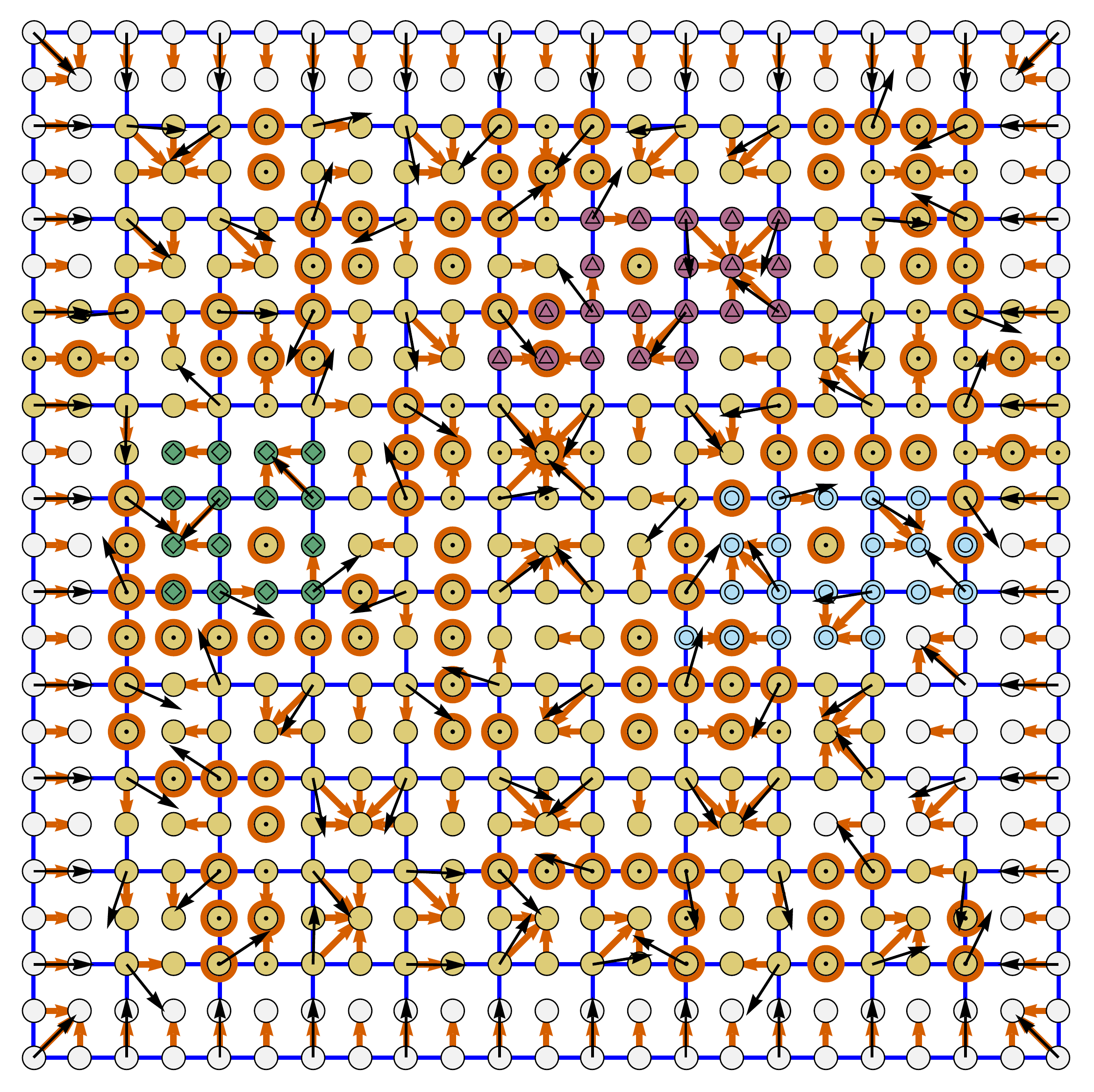}
\end{center}
  \caption{A combinatorial multivector field constructed from a random collection of vectors
  at the lattice points.}
  \label{fig:random-cmvf}
\end{figure}

\subsection{A combinatorial multivector field constructed from a random collection of vectors.}
Figure~\ref{fig:random-cmvf} presents a combinatorial multivector field constructed by algorithm CMVF from
a randomly selected collection of vectors at the lattice points. To ensure that the boundary of the selected region
does not divide multivectors, all the vectors at the boundary  are not random but point inwards.
The resulting Morse decomposition  consists of 102  isolated invariant sets out of which three consist of more than one multivector.

\section{Preliminaries}
\label{sec:preli}

In this section we introduce the notation, recall the definitions
and gather results used in the sequel.

\subsection{Sets and maps.}
We denote  the sets of reals, integers, non-negative
integers and non-positive integers respectively by $\RR, \ZZ$, $\ZZ^+$, $\ZZ^-$.
We also write $\ZZ^{\geq n}$, $\ZZ^{\leq n}$
respectively for integers greater or equal $n$ and less or equal $n$.
Given a set $X$,  we write $\card X$ for the number of elements of $X$ and
we denote by $\cP(X)$ the family of all subsets of $X$.
We write $f:X\pto Y$ for a partial map from $X$ to $Y$, that is a map
defined on a subset $\dom f\subset X$, called the {\em domain} of $f$,
and such that the set of values of $f$, denoted
$\im f$, is contained in $Y$.

\subsection{Relations, multivalued maps and digraphs.}
Given a set $X$ and a binary relation $R\subset X\times X$,
we use the shorthand $xRy$ for $(x,y)\in R$.
By the {\em transitive closure} of $R$ we mean
the relation $\bar{R}\subset X\times X$
given by $x\bar{R}y$ if there exists
a sequence $x=x_0,x_1,\ldots, x_n=y$ such that $n\geq 1$ and  $x_{i-1}Rx_i$
for $i=1,2,\ldots, n$. Note that $\bar{R}$ is transitive but need not be reflexive.
The  relation
$\bar{R}\cup\id_X$, where $\id_X$ stands for the identity relation on $X$,
is reflexive and transitive. Hence,
it is a {\em preorder}, called the {\em preorder induced by $R$}.
A $y\in X$ {\em covers} an $x\in X$ in the relation $R$ if $xRy$ but there is no $z\in X$
such that $x\neq z\neq y$ and $xRz$, $zRy$.

A multivalued map $F:X\mvmap Y$ is a map $F:X\to\cP(Y)$.
For $A\subset X$ we define the {\em image of $A$} by
$
    F(A):=\bigcup \setof{F(x) \mid x\in A}
$
and for $B\subset Y$ we define the {\em preimage of $B$} by
$
    F^{-1}(B):=\setof{x\in X\mid F(x)\cap B\neq\emptyset}.
$

Given a relation $R$, we associate with it a multivalued map $F_R:X\mvmap X$,
by $F_R(x):=R(x)$, where
$
   R(x):=\setof{y\in X\mid xRy}
$
is the {\em image of $x\in X$ in $R$}.
Obviously $R\mapsto F_R$  is a one-to-one correspondence between
binary relations in $X$ and multivalued maps from $X$ to $X$.
Often, it will be convenient to interpret the relation $R$ as a directed graph
whose set of vertices is $X$ and a directed arrow joins $x$ with $y$ whenever $xRy$.
The three concepts relation, multivalued map and directed graph are
equivalent on the formal level and the distinction is used only to emphasize different directions of research.
However, in this paper it will be convenient to use all these concepts interchangeably.

\subsection{Partial orders.}

Assume $(X,\leq)$ is a poset. Thus, $\leq$ is a {\em partial order}, that is
a reflexive, antisymmetric and transitive relation in $X$.
As usual, we denote the inverse of this relation by $\geq$. We also write $<$ and $>$ for the associated
strict partial orders, that is relations $\leq$ and $\geq$ excluding identity.
By an {\em interval} in $X$ we mean a subset of $X$ which has one of the following four forms
\begin{eqnarray*}
    ~~[x,y]&:=&\setof{z\in X\mid x\leq z\leq y}\\
    ~~(-\infty,y]&:=&\setof{z\in X\mid z\leq y}\\
    ~~[x,\infty)&:=&\setof{z\in X\mid x\leq z}\\
    ~~(-\infty,\infty)&:=&X.
\end{eqnarray*}
In the first case we speak about a {\em closed interval}. The elements $x,y$ are the {\em endpoints} of the interval.
We recall that $A\subset X$ is {\em convex} if for any $x,y\in A$ the closed interval $[x,y]$ is contained in $A$.
Note that every interval is convex but there may exists convex subsets of $X$ which are not intervals.
A set $A\subset X$ is an {\em upper set}
if for any $x\in X$ we have $[x,\infty)\subset A$. Also, $A\subset X$ is  a {\em lower set}
if for any $x\in X$ we have $(-\infty,x]\subset A$.
Sometimes a lower set is called an attracting interval and an upper set a repelling interval.
However, one has to be careful, because in general lower and upper sets need not be intervals at all.
For $A\subset X$ we also use the notation
  $A^{\leq}:=\setof{x\in X\mid \exists_{a\in A}\;x\leq a}$ and
  $A^{<}:=A^{\leq}\setminus A$.

\begin{prop}
\label{prop:interval}
   If $I$ is convex, then $I^{\leq}$ and $I^{<}$ are lower sets (attracting intervals).
\end{prop}
\proof
  The verification that  $I^{\leq}$ is a lower set is straightforward.
To see that $I^{<}$ is a lower set take $x\in I^{<}$. Hence, we have $x\not\in I$ but
$x<z$ for some $z\in I$. Let $y\leq x$. Then  $y\in I^{\leq}$.
Since $I$ is convex, we cannot have $y\in I$. It follows that $y\in I^{<}$.
\qed

\begin{prop}
\label{prop:lower-integer}
Let $I=\setof{1,2,\ldots n}$ and let $\leq$ denote the linear order of natural numbers.
The for any $i\in I$ we have $\{i\}^{\leq}=\setof{1,2,\ldots i}$ and
$\{i\}^{<}=\setof{1,2,\ldots i-1}$.
\qed
\end{prop}

\subsection{Topology of finite sets.}

For a topological space $X$ and $A\subset X$ we write $\cl A$ for the closure of $A$.
We also define the {\em mouth} of $A$ by
\[
        \mouth A:=\cl A\setminus A.
\]
Note that $A$ is closed if and only if its mouth is empty.
We say that $A$ is {\em proper} if $\mouth A$ is closed.
Note that open and closed subsets of $X$ are proper.
In the case of finite topological spaces proper sets have a special structure.
To explain it we first recall some properties of finite topological spaces
based on the following fundamental result which goes back to P.S.~Alexandroff~\cite{Al1937}.
\begin{thm}
\label{thm:alexandroff}
For a finite poset $(X,\leq)$ the family $\cT_{\leq}$ of upper sets of $\leq$
is a $T_0$ topology on $X$. For a finite $T_0$ topological space $(X,\cT)$ the relation
$x\leq_{\cT}y$ defined by $x\in\cl\{y\}$ is a partial order on $X$.
Moreover, the two associations relating $T_0$ topologies and partial orders
are mutually inverse.
\end{thm}

Let $(X,\cT)$ be a finite topological space. For $x\in X$ we write
 $ \cl x:=\cl\{x\}$,
  $\opn x:=\bigcap\setof{U\in\cT\mid x\in U}$.
The following proposition may be easily verified.
\begin{prop}
\label{prop:f-topology}
Let $(X,\cT)$ be a finite topological space and let $x,y\in X$.
The operations $\cl$ and $\opn$
have the following properties.
\begin{itemize}
   \item[(i)] $\cl x$ is the smallest closed set containing $x$,
   \item[(ii)] $\opn x$ is the smallest open set containing $x$,
   \item[(iii)] $\cl A=\bigcup_{x\in A}\cl x$ for any $A\subset X$,
   \item[(iv)] $A\subset X$ is closed if and only if
                $\cl x\subset A$ for any $x\in A$,
   \item[(v)] $A\subset X$ is open if and only if
                $\opn x\subset A$ for any $x\in A$,
   \item[(vi)] $y\in\cl x$  if and only if
                $x\in\opn y$.\qed
\end{itemize}
\end{prop}

In the sequel we will particularly often use property (iii) of Proposition~\ref{prop:f-topology}.
In particular, it is needed in the following characterization of proper sets in finite topological spaces.

\begin{prop}
\label{prop:proper-sets}
Let $X$ be a finite topological space. Then $A\subset X$ is proper if and only if
\begin{equation}
\label{eq:proper-sets}
\forall_{x,z\in A}\forall_{y\in X} \quad x\in \cl y, y\in \cl z\implies y\in A.
\end{equation}
\end{prop}
\proof
  Let $A\subset X$ be proper, $x,z\in A$, $y\in X$, $x\in \cl y$, $y\in \cl z$
and assume $y\not\in A$. Then $y\in \mouth A$ and $x\in\cl\mouth A=\mouth A$.
Therefore, $x\not\in A$, a contradiction proving \eqref{eq:proper-sets}.
Assume in turn that \eqref{eq:proper-sets} holds and $\mouth A$ is not closed.
Then  there exists an $x\in\cl\mouth A\setminus\mouth A$.
Thus, $x\in A$, $x\in\cl y$ for some $y\in\mouth A$ and $y\in\cl z$ for some $z\in A$.
It follows from \eqref{eq:proper-sets} that $y\in A$, which contradicts $y\in \mouth A$.
\qed

Proposition~\ref{prop:proper-sets} means that in the setting of finite topological spaces
proper sets correspond to convex sets in the language of the associated partial order.

\subsection{Graded modules and chain complexes}

Let $R$ be a fixed ring with unity.
Given a set $X$ we denote by $R(X)$ the free module over $R$
spanned by $X$.
Given a graded, finitely generated module $E=(E_k)_{k\in{\scriptsize \mathbb{Z^+}}}$ over $R$, we write
\[
     p_E(t):=\sum_{k=0}^{\infty} \rank(E_k)t^k,
\]
for the Poincar\'e formal power series of $E$.
We have the following theorem (see \cite{RZ1985})
\begin{thm}
\label{thm:exact-sequence}
Assume $E,F,G$ are graded, finitely generated modules and we have an exact sequence
\[
   \begin{diagram}
  \dgARROWLENGTH 5mm
    \node{\ldots}
    \arrow{e,t}{\gamma_{i+1}}
    \node{E_i}
    \arrow{e,t}{\alpha_{i}}
    \node{F_i}
    \arrow{e,t}{\beta_{i}}
    \node{G_i}
    \arrow{e,t}{\gamma_{i}}
    \node{\ldots}
    \arrow{e,t}{\gamma_{1}}
    \node{E_0}
    \arrow{e,t}{\alpha_{0}}
    \node{F_0}
    \arrow{e,t}{\beta_{0}}
    \node{G_0}
    \arrow{e,t}{\gamma_{0}}
    \node{0.}
   \end{diagram}
\]
Then
\begin{equation*}
p_E(t)+p_G(t)=p_F(t)+(1+t)Q(t),
\end{equation*}
where
\[
  Q(t):=\sum_{k=0}^{\infty} \rank(\im\gamma_{k+1})t^k
\]
is a polynomial with non-negative coefficients. Moreover, if $F=E\oplus G$, then $Q=0$.
\qed
\end{thm}

\section{Multivector fields and multivector dynamics.}
\label{sec:multi}

In this section we define Lefschetz complexes and
introduce the concepts of the combinatorial multivector and the combinatorial vector field on a Lefschetz complex.
Given a combinatorial multivector field,  we associate with it a graph and a multivalued map allowing us to study
its dynamics. We also prove a crucial theorem about acyclic combinatorial multivector fields.

\subsection{Lefschetz complexes}
\label{sec:lefschetz}
The following definition
goes back to S. Lefschetz
(see \cite[Chpt. III, Sec. 1, Def. 1.1]{Le1942}).

\begin{defn}
{\em
We say that $(X,\kappa)$ is a {\em Lefschetz complex}
if $X=(X_q)_{q\in{\scriptsize \mathbb{Z}^+}}$ is a finite set with gradation,
$\kappa : X \times X \to R$ is a map such that $\kappa(x,y)\neq 0 $
implies $x\in X_q,\,y\in X_{q-1}$ and
and for any $x,z\in X$ we have
\begin{equation}
\label{eq:kappa-condition}
    \sum_{y\in X}\kappa(x,y)\kappa(y,z)=0.
\end{equation}
We refer to the elements of $X$ as {\em cells}
and to $\kappa(x,y)$ as the {\em incidence coefficient} of $x,y$.
}
\end{defn}

The family of cells of a simplicial complex \cite[Definition 11.8]{KMW2004} and the family of
elementary cubes of a cubical set \cite[Definition 2.9]{KMW2004} provide simple but important examples of Lefschetz complexes.
In these two cases the respective formulas for the incident coefficients are explicit and elementary (see \cite{MB2009}).
In the case of a general regular cellular complex (regular finite CW complex, see \cite[Section IX.3]{Ma1991})
the incident coefficients may be obtained
from a system of equations (see \cite[Section~IX.5]{Ma1991}).

The Lefschetz complex $(X,\kappa)$ is called {\em regular} if for any $x,y\in X$
the incidence coefficient
$\kappa(x,y)$ is either zero or is invertible in $R$.
One easily verifies that condition \eqref{eq:kappa-condition}
implies that we have a free chain complex $(R(X),\bdy^\kappa)$ with
$\bdy^\kappa:R(X)\to R(X)$ defined
on generators by $\bdy^\kappa(x) := \sum_{y\in X}\kappa(x,y)y$.
The {\em Lefschetz homology} of $(X,\kappa)$,
denoted $H^\kappa(X)$, is the homology of this chain complex.
By a {\em zero space} we mean a Lefschetz complex whose Lefschetz homology is zero.
Since $X$ is finite, $(R(X),\bdy^\kappa)$ is finitely generated.
In consequence, the Poincar\'e formal power series $p_{H^{\kappa}(X)}(t)$
is a polynomial. We denote it briefly by $p_X(t)$.

Given $x,y\in X$ we say that $y$ is a {\em facet} of $x$ and write $y\adhl_{\kappa} x$ if $\kappa(x,y)\neq 0$.
It is easily seen that the relation $\adhl_{\kappa}$ extends uniquely to a minimal partial order.
We denote this partial order by $\leq_\kappa$
and the associated strict order by $<_\kappa$.
We say that $y$ is a {\em face} of $x$ if $y\leq_\kappa x$.
The $T_0$ topology defined via Theorem~\ref{thm:alexandroff} by the partial order
$\leq_\kappa$ will be called the
{\em Lefschetz topology} of $(X,\kappa)$. Observe that the closure of a set $A\subset X$
in this topology consists of all faces of all cells in $A$.
The Lefschetz complex via its Lefschetz topology may be viewed as an example of the abstract
cell complex in the sense of \cite{Ko1989}. It is also related to the abstract
cell complex in the sense of  \cite[Section III]{St1908}).

\begin{prop}
\label{prop:Lef-sing-dubl}
If $X=\{a\}$ is a singleton, then $H^\kappa(X)\cong R(X)\neq 0$.
If $X=\{a,b\}$ and $\kappa(b,a)$ is invertible, then $H^\kappa(X)=0$.
\end{prop}
\proof
If $X=\{a\}$, then $\bdy^{\kappa}$ is zero.
If $X=\{a,b\}$ and $\kappa(b,a)$ is invertible, then the only non-zero component of $\bdy^{\kappa}$ is an isomorphism.
\qed

Proposition~\ref{prop:Lef-sing-dubl} shows that a Lefschetz complex consisting
of just two cells may have zero Lefschetz homology. At the same time
the singular homology of this two point space with Lefschetz topology is non-zero,
because the singular homology of a non-empty space is never zero.
Thus, the singular homology $H(X)$ of a Lefschetz complex $(X,\kappa)$
considered as a topological space with its Lefschetz topology need
not be the same as the Lefschetz homology $H^\kappa(X)$.

A set $A\subset X$ is a {\em $\kappa$-subcomplex} of $X$
if $(A,\kappa_{|A\times A})$ is a Lefschetz complex.
Lefschetz complexes, under the name of S-complexes, are discussed in
\cite{MB2009}.
In particular, the following proposition follows from
the observation that a proper subset of a Lefschetz complex $X$
satisfies the assumptions of \cite[Theorem 3.1]{MB2009}.
\begin{prop}
\label{prop:kappa-subcomplex}
Every proper $A\subset X$  is a $\kappa$-subcomplex of $X$.
In particular, open  and  closed subsets of $X$ are $\kappa$-subcomplexes of $X$.
\qed
\end{prop}

Note that  a $\kappa$-subcomplex $A$ of $X$  does not guarantee
that $(R(A),\bdy^\kappa_{|R(A)})$ is a chain subcomplex of  $(R(X),\bdy^\kappa)$.
However, we have the following theorem (see \cite[Theorem 3.5]{MB2009}).
\begin{thm}
\label{thm:closed-subcomplex}
Assume $A$ is closed in $X$. Then $(R(A),\bdy^\kappa_{|R(A)})$ is a chain subcomplex of
$(R(X),\bdy^\kappa)$. Moreover, the homomorphisms $\bdy^{\kappa_{|A\times A}}:R(A)\to R(A)$
and $\bdy^\kappa_{|R(A)}:R(A)\to R(A)$ coincide. In particular, the  homology
of the quotient chain complex $(R(X)/R(A),[\bdy^\kappa])$, denoted
$H^\kappa(X,A)$, is well defined and isomorphic to $H^\kappa(X\setminus A)$.
\qed
\end{thm}

The following proposition is straightforward to verify.
\begin{prop}
\label{prop:hom-union}
Assume $X=X_1\cup X_2$, where $X_1$ and $X_2$ are disjoint, closed subset of $X$.
Then $X_1$ and $X_2$ are $\kappa$-subcomplexes and $H^\kappa(X)=H^\kappa(X_1)\oplus H^\kappa(X_2)$.
\qed
\end{prop}

We will also need the following theorem which follows from \cite[Theorems 3.3 and 3.4]{MB2009}
\begin{thm}
\label{thm:exact-sequence-L}
   Assume $X'\subset X$ is closed in $X$, $\kappa':=\kappa_{|X'\times X'}$ and $X'':=X\setminus X'$.
Then  there is a long exact sequence of homology modules
\begin{equation}
\label{eq:long-exact-sequence}
  \ldots \mapright{ } H^\kappa_i(X') \mapright{ } H^\kappa_i(X) \mapright{ }
  H^\kappa_i(X'') \mapright{ } H^\kappa_{i-1}(X') \mapright{ } \ldots ~~.
\end{equation}
\end{thm}

\subsection{Multivectors.}
Let $(X,\kappa)$ be a fixed Lefschetz complex.
\begin{defn}{\em
A {\em combinatorial multivector} or briefly a {\em multivector} is a proper subset $V\subset X$ admitting
a unique maximal element with respect to the partial order $\leq_\kappa$.
We call this element the  {\em dominant} cell of $V$ and denote it $V^\star$.
}
\end{defn}
Note that we do not require the existence of a unique minimal element in a multivector
but if such an element exists, we denote it by $V_\star$.
Multivectors admitting a unique minimal element are studied in \cite{Freij2009}
in the context of equivariant discrete Morse theory.
A concept similar to our multivector appears also in \cite{Wisn2008}.

\begin{prop}
\label{prop:multivector}
For a multivector $V$ we have
  $V\neq\emptyset$ and  $\cl V=\cl V^\star$.\qed
\end{prop}

A multivector is {\em regular} if $V$ is a zero space.
Otherwise it is called {\em critical}.
A combinatorial multivector $V$  is a {\em combinatorial vector} or briefly a {\em vector}  if $\card V\leq 2$.
A vector always has a unique minimal element.

\begin{prop}
\label{prop:vector}
Assume $X$ is a regular Lefschetz complex and let
$V\subset X$ be a vector. Then $\card V=1$ if and only if $V$ is critical
and $\card V=2$ if and only if $V$ is regular. Moreover, if $\card V=2$,
then $V_\star\adhl_{\kappa} V^\star$.
\end{prop}
\proof
   If $V$ is a singleton, then by Proposition~\ref{prop:Lef-sing-dubl} we have
$H^\kappa(V)\neq 0$, hence $V$ is critical.
If $\card V=2$, then $V_\star \neq V^\star$.
First, we will show that $V_\star\adhl_{\kappa} V^\star$.
Indeed, if not, then $V_\star <_{\kappa} x <_{\kappa} V^\star$ for some $x\in X$. But then $x\in\mouth V$
and $V_\star\in \cl\mouth V\setminus\mouth V$, which contradicts the assumption that $V$ is proper.
Thus, $V_\star$ is a facet of $V^\star$, $\kappa(V^\star,V_\star)\neq 0$ and by the assumed regularity of $X$
it is invertible. Therefore, again by Proposition~\ref{prop:Lef-sing-dubl}, we have
$H^\kappa(V)=0$. It follows that $V$ is regular.
\qed

\subsection{Multivector fields.}
The following definition is the main new concept introduced in this paper.
\begin{defn}{\em
A {\em combinatorial multivector field} on $X$, or briefly a {\em multivector field},
is a partition $\cV$ of $X$ into multivectors.
A {\em combinatorial vector field} on $X$, or briefly a {\em vector field},
is a combinatorial multivector field whose each multivector is a vector.
}
\end{defn}
Proposition~\ref{prop:vector} implies that our concept of a vector field on the Lefschetz complex
of a cellular complex is in one-to-one correspondence
with Forman's combinatorial vector field (see \cite{Fo98b}).
It also corresponds to the concept of partial matching \cite[Definition 11.22]{Ko}.
Thus, the combinatorial multivector field is a generalization of the earlier definitions
in which vectors were used instead of multivectors.

For each cell $x\in X$ we denote by $\vclass{x}$
the unique multivector in $\cV$ to which $x$ belongs.
If the multivector field $\cV$ is clear from the context,
we write briefly $[x]:=\vclass{x}$ and $x^\star:=\vclass{x}^\star$.
We refer to a cell $x$ as {\em dominant with respect to $\cV$}, or briefly as {\em dominant},
if $x^\star=x$.

The map which sends $x$ to $x^\star$ determines the combinatorial multivector field.
More precisely, we have the following theorem.

\begin{thm}
\label{thm:theta}
The map $\theta:X\ni x\mapsto x^\star\in X$ has the following properties
\begin{itemize}
   \item[(i)] for each $x\in X$ we have $x\in\cl \theta(x)$,
   \item[(ii)] $\theta^2=\theta$,
   \item[(iii)] for each $y\in\im\theta$ if $x\in \theta^{-1}(y)$, then $\opn x\cap\cl y\subset \theta^{-1}(y)$.
\end{itemize}
Conversely, if a map $\theta:X\to X$ satisfies properties (i)-(iii), then
\[
\cV_\theta:=\setof{\theta^{-1}(y)\mid y\in\im\theta}
\]
is a combinatorial multivector field on $X$.
\end{thm}
\proof
Properties (i)-(iii) of $\theta:X\ni x\mapsto x^\star\in X$
follow immediately from the definition of a multivector.
To prove the converse assertion assume $\theta:X\to X$ satisfies properties (i)-(iii).
Obviously $\cV_\theta$ is a partition. To see that each element of $\cV_\theta$
is a multivector take $y\in \im\theta$. Then by (i) $\theta^{-1}(y)\subset \cl y$
and by (iii) $\theta^{-1}(y)$ is open in $\cl y$.
This proves that $\theta^{-1}(y)$ is proper.
By (ii) the unique maximal element in $\theta^{-1}(y)$ is $y$.
Therefore, $\theta^{-1}(y)$ is a multivector.
\qed

\subsection{The graph and multivalued map of a multivector field.}
Given a combinatorial multivector field $\cV$ on $X$ we associate with it
the graph $G_\cV$ with vertices in $X$ and an arrow from $x$ to $y$ if one of the following
conditions is satisfied
\begin{eqnarray}
  && x\neq y=x^\star \text{ (an {\em up-arrow}),}\label{eq:adhr1}\\
  && x=x^\star \text{ and } y\in\cl x\setminus \vclass{x}  \text{ (a {\em down-arrow}),}\label{eq:adhr2}\\
  && x=x^\star=y \text{ and } [y] \text{ is critical (a {\em loop}).}\label{eq:adhr3}
\end{eqnarray}
We write $y \adhl_{\cV} x$ if there is an arrow from $x$ to $y$ in $G_\cV$.
This lets us interpret $\adhl_{\cV}$ as a relation in $X$.
We denote by $\leqV$ the preorder induced by $\adhl_{\cV}$.
In order to study the dynamics of $\cV$, we interpret $\adhl_{\cV}$
as a multivalued map $\Pi_{\cV}: X\mvmap X$, which
sends a cell $x$ to the set of cells covered by $x$ in $\adhl_{\cV}$, that is
\[
  \Pi_\cV(x):=\setof{y\in X\mid y\adhl_{\cV} x}.
\]

We say that a cell $x\in X$ is {\em critical} with respect to $\cV$
if $\vclass{x}$ is critical and $x$ is dominant in $\vclass{x}$.
A cell is {\em regular} if it is not critical. We denote by $\regclass{x}$ the set of regular cells
in $[x]_\cV$. It is straightforward to observe that
\begin{equation}
\label{eq:regclass}
   \regclass{x}=\begin{cases}
                  \vclass{x} & \text{if $\vclass{x}$ is regular,}\\
                  \vclass{x}\setminus\{x^\star\} & \text{otherwise.}
                \end{cases}
\end{equation}

\begin{prop}
\label{prop:Pi_cV}
  For each $x\in X$ we have
\[
   \Pi_\cV(x)=\begin{cases}
         \{x^\star\} & \text{if $x\neq x^\star,$}\\
         \cl x\setminus\regclass{x} & \text{otherwise.}\qed
             \end{cases}
\]
\end{prop}

We extend the relation $\leqV$ to multivectors $V,W\in\cV$ by assuming that
$V\leqV W$ if and only if  $V^\star\leqV W^\star$.

\begin{prop}
\label{prop:cV-extension}
If $\leqV$ is a partial order on $X$, then the extension of $\leqV$ to multivectors
is a partial order on $\cV$.
\qed
\end{prop}

\subsection{Acyclic multivector fields.}
We say that $\cV$ is {\em acyclic} if $\leqV$
is a partial order on $X$.

\begin{thm}
\label{thm:acyclic-gdvf}
   Assume $X$ admits an acyclic multivector field whose each multivector is regular.
   Then $X$ is a zero space.
\end{thm}
\proof
  Let $\cV$ be an acyclic multivector field on $X$ whose each multivector is regular.
We will proceed by induction on $\card\cV$. If $\card \cV=0$, that is if $X$ is empty,
the conclusion is obvious. Assume $\card\cV>0$.
By Proposition~\ref{prop:cV-extension} we know that $\leqV$ is a partial order on $\cV$.
Let $V$ be a maximal element of $\cV$ with respect to $\leqV$.
We claim that $X':=X\setminus V$ is closed in $X$. To prove the claim, assume the contrary.
Then  there exists an $x\in\cl X'\cap V$.
Let $y\in X'$ be such that $x\in\cl y\subset \cl y^\star$. Since $x\in V$ and $y\not\in V$,
we have $\vclass{x}\neq\vclass{y}=\vclass{y^\star}$. It follows that $x\leqV y^\star$ and consequently
$V=\vclass{x}\leqV \vclass{y^\star}$. Hence, $V$ is not maximal, because $V\neq\vclass{y^\star}$, a contradiction proving that $X'$ is closed. In particular $X'$ is proper.
Obviously, $\cV':=\cV\setminus\{V\}$ is
an acyclic  multivector field on $X'$ whose each multivector is regular.
Thus, by induction assumption, $X'$ is a zero space.
Since also $V$ is a zero space, it follows from Theorem~\ref{thm:exact-sequence-L}
applied to the pair $(X,X')$
that $X$ is a zero space.
\qed

\section{Solutions and invariant sets}
\label{sec:solu}

In this section we first define the solution of a combinatorial multivector field,
an analogue of a solution of an ordinary differential equation.
Then, we use it to define the fundamental concept of the invariant set
of a combinatorial multivector field.

\subsection{Solutions.}
A partial map $\varphi:\ZZ\pto X$ is a {\em solution} of $ \cV$ if $\dom\varphi$
is an interval in $\ZZ$ and
$
\varphi(i+1)\in\Pi_\cV(\varphi(i)) \text{ for } i,i+1\in\dom\varphi.
$
A solution $\varphi$ is in $A\subset X$ if $\im\varphi\subset A$.
We call $\varphi$ a {\em full} (respectively {\em forward} or {\em backward}) solution
if $\dom\varphi$ is $\ZZ$ (respectively $\ZZ^{\geq n}$ or $\ZZ^{\leq n}$ for some $n\in\ZZ$).
We say that $\varphi$ is a {\em solution through $x\in X$} if $x\in\im\varphi$.
We denote the set of full (respectively forward, backward)
solutions in $A$ through $x$ by $\Sol(x,A)$
(respectively $\Sol^+(x,A)$, $\Sol^-(x,A)$).
We drop $A$ in this notation if $A$ is the whole Lefschetz complex $X$.
As an immediate consequence of Proposition~\ref{prop:Pi_cV} we get the following proposition.
\begin{prop}
\label{prop:solution-star}
If $\varphi$ is a solution of $\cV$ and $i,i+1\in\dom\varphi$ then either
$\varphi(i)^\star=\varphi(i)$ or $\varphi(i)^\star=\varphi(i+1)$.
\qed
\end{prop}

Given $n\in\ZZ$, let
$
  \tau_n:\ZZ\ni i\mapsto i+n\in\ZZ
$
denote the $n$-translation map.
Let $\varphi$ be a solution in $A$ such that $n\in\ZZ$ is the right endpoint of $\dom\varphi$
and let $\psi$ be a solution  in $A$ such that $m\in\ZZ$ is the left endpoint of $\dom\psi$.
We define
$
\varphi\cdot\psi:\tau_n^{-1}(\dom\varphi)\cup \tau_{m-1}^{-1}(\dom\psi)\to A,
$
the {\em concatenation} of $\varphi$ and $\psi$ by
\[
  (\varphi\cdot\psi)(i):=
     \begin{cases}
     \varphi(i+n)  &\text{if $i\leq 0$}\\
     \psi(i+m-1) & \text{if $i>0$}.
     \end{cases}
\]
It is straightforward to observe that if  $\psi(m)\in\Pi_\cV(\varphi(n))$
then the concatenation $\varphi\cdot\psi$ is also a solution in $A$.

Let $\varphi\in\Sol^+(x,A)$ and let $\dom\varphi=\ZZ^{\geq n}$ for some $n\in\ZZ$.
We define $\sigma_+\varphi:\ZZ^{\geq n}\to A$, the {\em right shift} of $\varphi$,
by $\sigma_+\varphi(i):=\varphi(i+1)$.
Let $\psi\in\Sol^-(x,A)$ and let $\dom\psi=\ZZ^{\leq n}$ for some $n\in\ZZ$.
We define $\sigma_-\psi:\ZZ^{\leq n}\to A$, the {\em left shift} of $\psi$, by
$\sigma_-\psi(i):=\psi(i-1)$.
It is easily seen that
for $k\in\ZZ^+$ we have $\sigma_+^k\varphi\in\Sol^+(\varphi(n+k),A)$
and $\sigma_-^k\psi\in\Sol^-(\psi(n-k),A)$.
For a full solution $\varphi$ we denote respectively by $\varphi^+$, $\varphi^-$
the restrictions $\varphi_{|{\scriptsize \mathbb{Z}^+}}$,  $\varphi_{|{\scriptsize \mathbb{Z}^-}}$.

Obviously, if $\varphi$ is a solution in $A$
then $\varphi\circ\tau_n$ is also a solution in $A$.
We say that  solutions $\varphi$ and $\varphi'$ are {\em equivalent}
if $\varphi'=\varphi\circ\tau_n$ or $\varphi=\varphi'\circ\tau_n$ for some $n\in\ZZ$.
It is straightforward to verify that this is indeed an equivalence relation.
It preserves forward, backward and full solutions through $x$.
Moreover, it is not difficult to verify that the concatenation extends
to an associative operation on equivalence classes of solutions.
In the sequel we identify solutions in the same equivalence class.
This allows us to treat the solutions as finite or infinite words over the alphabet
consisting of cells in $X$.
In the sequel, whenever we pick up a representative of a forward (backward) solution, we assume
its domain is respectively $\ZZ^+$, ($\ZZ^-$).

\subsection{Paths.}
A solution $\varphi$ such that $\dom\varphi$ is a finite interval
is called a {\em path} joining the value of $\varphi$ at the left end of the domain
with the value at the right end. The cardinality of $\dom\varphi$ is called
the {\em length} of the path.
In the special case when $\varphi$ has
length one we identify $\varphi$ with its unique value.
We also admit the trivial path of length zero (empty set).
In particular, it acts as the neutral element of concatenation.
For every $x\in X$ there is a unique path  joining $x$ with $x^\star$,
denoted $\nu(x)$ and given by
\[
  \nu(x):=\begin{cases}
               x & \text{if $x=x^\star$,}\\
               x\cdot x^\star & \text{otherwise.}
             \end{cases}
\]
Note that the concatenation $x\cdot x^\star$ is a solution, because
$x^\star\in\Pi_\cV(x)$.
We also define
\[
  \nu^-(x):=\begin{cases}
               \emptyset & \text{if $x=x^\star$,}\\
               x & \text{otherwise.}
             \end{cases}
\]
Note that $\nu(x)=\nu^-(x)\cdot x^\star$.

\subsection{$\cV$-compatibility.}
We say that $A\subset X$ is $ \cV$-{\em compatible} if
$x\in A$ implies  $\vclass{x}\subset A$   for   $x\in X$.
We denote by $[A]_\cV^-$ the maximal $\cV$-compatible subset of $A$ and by
$[A]_\cV^+$ the minimal $\cV$-compatible superset of $A$.
The following proposition is straightforward.
\begin{prop}
\label{prop:V-compatibility}
For any $A,B\subset X$ we have
\begin{itemize}
   \item[(i)] $[A]_\cV^- \subset A \subset [A]_\cV^+$,
   \item[(ii)] if $A,B$ are $\cV$-compatible, then also $A\cap B$ and $A\cup B$ are $\cV$-compatible,
   \item[(iii)] $[A]_\cV^-$ and $[A]_\cV^+$ are $\cV$-compatible,
   \item[(iv)] $A\subset B$ implies $[A]_\cV^\pm\subset [B]_\cV^\pm$. \qed
\end{itemize}
\end{prop}
Observe that if $A\subset X$ is a $\cV$-compatible $\kappa$-subcomplex of $X$ then $\cV':=\setof{V\in\cV\mid V\subset A}$
is a multivector field on $A$. We call it the {\em restriction} of $\cV$ to $A$ and denote it $\cV_{|A}$.

\begin{figure}
\begin{center}
  \includegraphics[width=0.35\textwidth]{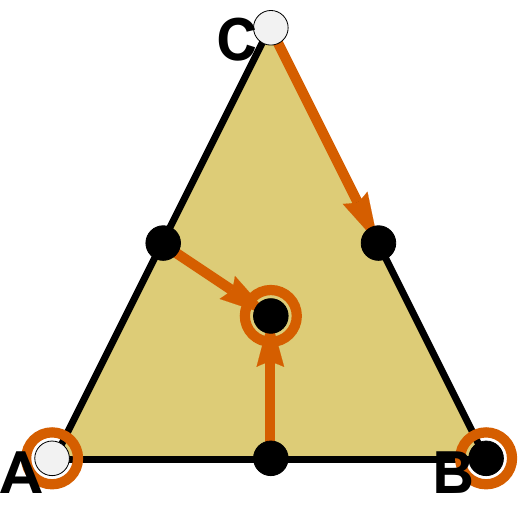}
\end{center}
  \caption{A multivector field on a triangle.
  The invariant part of the collection of cells $\{AB,AC,ABC,B,BC\}$ marked in black
  is $\{AB,AC,ABC,B\}$.}
  \label{fig:invSetproblems}
\end{figure}

\subsection{Invariant parts and invariant sets.}
We define the {\em invariant part of $A$}, the {\em positive invariant part of $A$}
and the {\em negative invariant part of $A$} respectively by
\begin{eqnarray*}
  \Inv A&:=& \setof{x\in A\mid \Sol(x^\star,[A]_\cV^-)\neq\emptyset},\\
  \Inv^+ A&:=& \setof{x\in A\mid   \Sol^+(x^\star,[A]_\cV^-)\neq\emptyset},\\
  \Inv^- A&:=& \setof{x\in A\mid \Sol^-(x^\star,[A]_\cV^-)\neq\emptyset}.
\end{eqnarray*}

Note that by replacing $x^\star$ by $x$ or $[A]_\cV^-$ by $A$ in the definition
of the invariant part of $A$ we may not obtain the invariant part of $A$.
Indeed, consider the set $\{AB,AC,ABC,B,BC\}$ marked in black in Figure~\ref{fig:invSetproblems}.
Its invariant part is $\{AB,AC,ABC,B\}$. But, by replacing $x^\star$ by $x$ in the definition
of the invariant part we obtain $\{ABC,B\}$. And by replacing $[A]_\cV^-$ by $A$
we obtain $\{AB,AC,ABC,B,BC\}$.

\begin{prop}
\label{prop:inv}
For any $A\subset X$ we have
\begin{eqnarray}
  [\Inv A]_{\cV}^-&=&\Inv A,\label{eq:inv-comp}\\
  \Inv A&\subset&A,\label{eq:inv0}\\
  A\subset B&\implies& \Inv A\subset \Inv B,\label{eq:inv0m}\\
  \Inv A&=&\Inv^- A \cap \Inv^+ A. \label{eq:inv1}\\
  \Inv A&=&\Inv \Inv A.  \label{eq:inv2}
\end{eqnarray}
\end{prop}
\proof
Equations~\eqref{eq:inv-comp},~\eqref{eq:inv0},~\eqref{eq:inv0m},~\eqref{eq:inv1}
and the right-to-left inclusion in  \eqref{eq:inv2}
are straightforward. To prove the left-to-right inclusion in \eqref{eq:inv2}
take $x\in\Inv A$ and let $\varphi\in\Sol(x^\star,[\Inv A]_{\cV}^-)$.
Fix an $i\in\ZZ$. Obviously,
\[
\varphi\in\Sol(\varphi(i),[A]_{\cV}^-)\cap\Sol(\varphi(i+1),[A]_{\cV}^-).
\]
Hence, by Proposition~\ref{prop:solution-star} $\varphi\in\Sol(\varphi(i)^\star,[A]_{\cV}^-)$,
which means that $\varphi(i)\in\Inv A$. It follows from~\eqref{eq:inv-comp} that $[\varphi(i)]_\cV\subset \Inv A$.
Thus, $\varphi(i)\in [\Inv A]_{\cV}^-$ and since $i\in\ZZ$ is arbitrarily fixed we conclude that
$\varphi\in\Sol(x^\star,[\Inv A]_{\cV}^-)\neq\emptyset$.
Therefore, $x\in\Inv \Inv A$.
\qed

We say that $A$ is  {\em invariant with respect to $\cV$} if $\Inv A=A$.
\begin{prop}
\label{prop:inv-set}
Assume $A\subset X$ is invariant.
Then $A$ is $\cV$-compatible. Moreover, for any $x\in A$
we have $\Sol^+(x,[A]_{\cV}^-)\neq\emptyset$ and for any dominant $x\in A$
we have $\Sol^-(x,[A]_{\cV}^-)\neq\emptyset$.
\end{prop}
\proof
   It follows from~\eqref{eq:inv-comp} that the invariant part of any set
is $\cV$-compatible. In particular, an invariant set, as the invariant part of itself
is $\cV$-compatible. Let $x\in A$ and $\varphi\in\Sol(x^\star,[A]_{\cV}^-)$.
Then $\nu^-(x)\cdot\varphi^+\in\Sol^+(x,[A]_{\cV}^-)$.
Obviously, $\varphi^-\in\Sol^-(x,[A]_{\cV}^-)$ when $x=x^\star$.
\qed

\section{Isolated invariant sets and the Conley index}
\label{sec:iso}

In this section we introduce the concept of an isolated invariant set
of a combinatorial multivector field and its homological invariant,
the Conley index. Both are analogues of the classical concepts for flows \cite{Co78}.
From now on we fix a combinatorial multivector field $\cV$ on a Lefschetz complex $X$
and we assume that $X$ is invariant with respect to $\cV$.

\subsection{Isolated invariant sets.}
Assume $A\subset X$ is invariant. We say that a path $\varphi$ from $x\in A$ to $y\in A$
is an {\em internal tangency} of $A$, if $\im \varphi\subset\cl A$ and $\im\varphi\cap\mouth A\neq\emptyset$.
We say that $S\subset X$ is an {\em isolated invariant set} if it is invariant and admits no internal tangencies.

\begin{thm}
\label{thm:iso-inv-set}
Let $S\subset X$ be invariant. Then $S$ is an isolated invariant set if and only if $S$ is proper.
\end{thm}
\proof
Assume $S\subset X$ is an isolated invariant set. By Proposition~\ref{prop:inv-set}
it is $\cV$-compatible. Assume to the contrary that $S$ is not proper.
Then there exists an $x\in\cl\mouth S\setminus\mouth S$.
Hence, $x\in\cl z$ for a $z\in\mouth S$ and $z\in\cl y$ for a $y\in S$.
In particular, $x\in S$ and $z\not\in S$. It follows from the $\cV$-compatibility of $S$
that $\class{x}\neq\class{z}\neq\class{y}$. Hence, $x\in\Pi_{\cV}(z)$ and $z\in\Pi_{\cV}(y)$.
Thus, $y\cdot z\cdot x$ is an internal tangency, a contradiction proving that $S$ is proper.

To prove the opposite implication, take $S\subset X$ which is invariant and proper
and assume to the contrary that $\varphi$ is an internal tangency of $S$.
Then the values $x,y$ of $\varphi$ at the endpoints of $\dom \varphi$ belong to $S$ and there is a $k\in\dom\varphi$
such that $\varphi(k)\in \mouth S$. Thus, we can choose an $m\in\dom\varphi$ satisfying
$\varphi(m)\not\in S$ and $\varphi(m+1)\in S$.
In particular, $\varphi(m)\in\mouth S$ and $\varphi(m+1)\not\in\mouth S$.
Proposition~\ref{prop:solution-star}
and $\cV$-compatibility of $S$ imply that $\varphi(m)^\star=\varphi(m)$.
It follows that $\varphi(m+1)\in\cl \varphi(m)\subset\cl\mouth S$.
However, this is  a contradiction, because $S$ is proper.
\qed

\subsection{Index pairs.}
A pair $P=(P_1, P_2)$ of closed subsets of $X$ such that $P_2\subset P_1$
is an {\em index pair} for $S$ if the following three conditions are satisfied
\begin{gather}
  P_1\cap\Pi_{\cV}(P_2)\subset P_2,\label{eq:ip1}\\
  P_1\cap\Pi_{\cV}^{-1}(X\setminus P_1)\subset P_2,\label{eq:ip2}\\
  S=\Inv( P_1\setminus P_2).\label{eq:ip3}
\end{gather}
We say that the index pair $P$ is {\em saturated} if  $P_1\setminus P_2=S$.

\begin{figure}
  \includegraphics[width=0.9\textwidth]{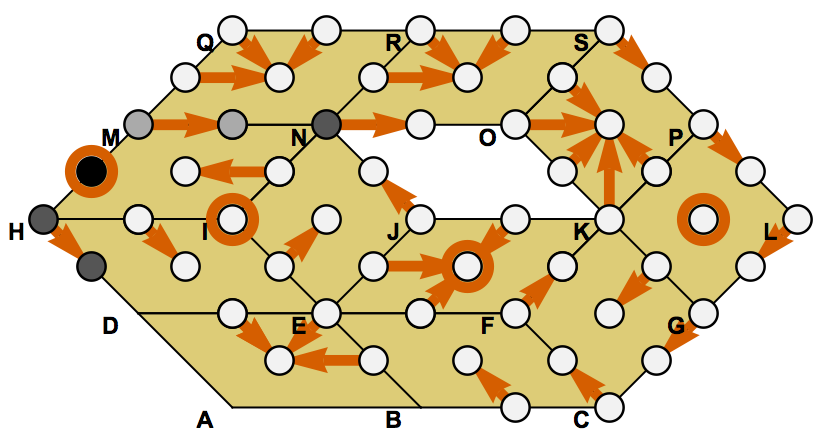}\\
  \caption{ The singleton consisting of the cell marked in black is a simple
  example of an isolated invariant set.
  Taking all six cells marked in non-white as $P_1$ and all three cells marked in dark gray as $P_2$
  we obtain a sample index pair $P=(P_1,P_2)$ of this isolated invariant set. This is not
  a saturated index pair. A saturated index pair may be obtained by adding to $P_2$ the two cells marked
  in light gray.}
  \label{fig:papex-iis2}
\end{figure}

A sample isolated invariant set together with an index pair is presented in Figure~\ref{fig:papex-iis2}.

\begin{prop}
\label{prop:ip-x-xstar}
If $(P_1,P_2)$ is an index pair and $x\in P_1\setminus P_2$ then $x^\star\in P_1\setminus P_2$.
\end{prop}
\proof
Assume to the contrary that $x^\star\not\in P_1\setminus P_2$.
Then  either $x^\star\not\in P_1$ or $x^\star\in P_2$.
Since $x^\star\in\Pi_{\cV}(x)$, the first case contradicts \eqref{eq:ip1}.
Since $x\in\cl x^\star$, the second case contradicts the fact that $P_2$ is closed.
\qed

\begin{prop}
\label{prop:ip-V-compatibility}
For any index pair $(P_1,P_2)$ the set $P_1\setminus P_2$ is $\cV$-compatible.
\end{prop}
\proof
Assume to the contrary that $x\in P_1\setminus P_2$ but $[x]\not\subset P_1\setminus P_2$.
By Proposition~\ref{prop:ip-x-xstar} we may assume without loss of generality that $x=x^\star$.
Then $y\not\in P_1\setminus P_2$ for some $y\in[x^\star]$.
Since $y\in\cl x^\star\subset P_1$, we see that $y\in P_2$. But $x^\star\in \Pi_{\cV}(y)$,
therefore we get from \eqref{eq:ip1} that $x^\star\in P_2$, a contradiction.
\qed\vspace*{1mm}

The following proposition follows immediately from the definition of index pair.
\begin{prop}
\label{prop:y-in-P1}
Assume $P=(P_1,P_2)$ is an index pair, $x\in P_1$ and $\varphi\in\Sol^+(x)$.
Then  either $\varphi\in\Sol^+(x,P_1)$ or $\varphi(i)\in P_2$ for some $i\in\ZZ^+$.
\qed
\end{prop}

Given an index pair $P$, consider the set
\[
  \E{P}:=\setof{x\in P_1\mid\forall_{ \varphi\in\Sol^+(x)}\,  \exists_{ i\in\dom\varphi}\;\varphi(i)\in P_2}
\]
of cells in $P_1$ whose every forward solution intersects $P_2$.

\begin{prop}
\label{prop:y-in-ep}
If $y\in \E{P}\setminus P_2$, then $y^\star\in \E{P}$.
\end{prop}
\proof
   If $y=y^\star$, then the conclusion is obvious. Therefore, we may assume that $y\neq y^\star$.
Then  $y^\star=\Pi_{\cV}(y)$ and since  $y\not\in P_2$ we get from \eqref{eq:ip2}
that $y^\star\in P_1$. Let $\varphi\in\Sol^+(y^\star)$.
Then  $\varphi':=y\cdot\varphi\in\Sol^+(y)$. The assumption $y\in \E{P}$ implies that
$\varphi'(i)\in P_2$ for some $i\in\ZZ^+$. Since $y\not\in P_2$ we have $i>0$.
It follows that $\varphi(i-1)=\varphi'(i)\in P_2$ and consequently $y^\star\in \E{P}$.
\qed

In the following lemma, given an index pair $P$ we first shrink $P_1$ and then grow $P_2$ to saturate $P$.
\begin{lem}
\label{lem:saturated}
  If $ P$ is an index pair for an isolated invariant set $S$, then
$
    P^*:=(S\cup \E{P}, P_2)
$
  is an index pair for $S$ and
$
    P^{**}:=(S\cup \E{P}, \E{P})
$
  is a saturated index pair for $S$.
\end{lem}
\proof
  First, we will show that the sets $\E{P}$ and $S\cup \E{P}$ are closed.
Take $x\in \cl \E{P}$ and assume $x\not\in\E{P}$.
Since $\E{P}\subset P_1$ and $P_1$ is closed, we have $x\in P_1$.
Choose $y\in \E{P}$ such that $x\in\cl y$.
If $y\in P_2$, then $x\in P_2\subset \E{P}$.
Hence, assume that $y\not\in P_2$.
Then, by Proposition~\ref{prop:y-in-ep}, $y^\star\in \E{P}$.
Moreover, $x\in\cl y\subset \cl y^\star$.
Since  $x\not\in \E{P}$,
there exists a $\varphi\in\Sol^+(x,X\setminus P_2)$.
But $x\in P_1$, hence it follows from Proposition~\ref{prop:y-in-P1}
that $\varphi\in\Sol^+(x,P_1\setminus P_2)$.
Note that $x\neq y^\star$, because $x\in\E{P}$ and $x\leq_\kappa y\leq_\kappa y^\star$.
Thus, we cannot have
$[x]=[y^\star]$, because otherwise $\varphi(1)=y^\star$ and
$\sigma(\varphi)\in\Sol^+(y^\star,P_1\setminus P_2)$, which contradicts $y^\star\in \E{P}$.
Hence, we have $[x]\neq [y^\star]$.
Then $x\in\Pi_{\cV}(y^\star)$ and $y^\star\cdot\varphi\in\Sol^+(y^\star,P_1\setminus P_2)$ which
again contradicts $y^\star\in \E{P}$. It follows that $x\in \E{P}$.
Therefore, $\E{P}$ is closed.

To show that $S\cup \E{P}$ is closed it is enough to prove that
\begin{equation}
\label{eq:cls-s}
\cl S\setminus S\subset \E{P}.
\end{equation}
Assume the contrary. Then  there exists
an $x\in \cl S\setminus S$ such that $x\not\in \E{P}$. Let $y\in S$ be such that $x\in\cl y$.
Without loss of generality we may assume that $y=y^\star$.
By the $\cV$-consistency of $S$ we have  $x\not\in[y^\star]$.
Thus, $x\in\Pi_{\cV}(y^\star)$.
Since $y^\star\in S\subset \Inv^-(P_1\setminus P_2)$,
we can take $\varphi\in\Sol^-(y^\star,P_1\setminus P_2)$.
Since $x\not\in \E{P}$
and $x\in\cl S\subset P_1$, we may take $\psi\in\Sol^+(x,P_1\setminus P_2)$.
Then  $\varphi\cdot\psi\in\Sol(x,P_1\setminus P_2)$.
By Proposition~\ref{prop:ip-V-compatibility} the set $P_1\setminus P_2$ is $\cV$-compatible.
Therefore,
$x\in \Inv(P_1\setminus P_2)=S$,
a contradiction again. It follows that $S\cup \E{P}$ is also closed.

In turn, we will show  that $P^*$ satisfies properties \eqref{eq:ip1}--\eqref{eq:ip3}.
Property \eqref{eq:ip1} follows immediately from the same property of index pair $P$, because
$S\cup \E{P}\subset P_1$.  To see \eqref{eq:ip2} take $x\in S\cup \E{P}$
and assume there exists a $y\in\Pi_{\cV}(x)\setminus (S\cup \E{P})$.
We will show first that $x\not\in S$. Assume to the contrary that $x\in S$.
It cannot be $y=x^\star$, because then $y\in S$
by the $\cV$-compatibility of $S$. Hence, $y\in\cl x$.
It follows that  $y\in \cl S\setminus S$ and by \eqref{eq:cls-s} $y\in \E{P}$,
a contradiction which proves that $x\not\in S$ and consequently $x\in \E{P}$.
Since $y\not\in\E{P}$, we can choose a solution $\varphi\in\Sol^+(y,X\setminus P_2)$.
It follows that $x\cdot\varphi\in\Sol^+(x)$. But, $x\in\E{P}$, therefore there exists
an $i\in\ZZ^+$ such that $(x\cdot\varphi)(i)\in P_2$. It must be $i=0$, because
for $i>0$ we have  $(x\cdot\varphi)(i)=\varphi(i-1)\not\in P_2$. Thus, $x=(x\cdot\varphi)(0)\in P_2$,
which proves \eqref{eq:ip2}.

Now, since $S\subset S\cup \E{P}\subset P_1$ and $S\cap P_2=\emptyset$, we have
\[
  S=\Inv S\subset \Inv(S\cup \E{P}\setminus P_2)\subset\Inv(P_1\setminus P_2)=S.
\]
This proves \eqref{eq:ip3} for $P^*$.

We still need to prove that $P^{**}$
satisfies properties \eqref{eq:ip1}-\eqref{eq:ip3}.
Let $x\in \E{P}$ and $y\in \Pi_{\cV}(x)\cap (S\cup \E{P})$. Then $y\in P_1$.
Let $\varphi\in \Sol^+(y)$. Then $x\cdot\varphi\in\Sol^+(x)$.
Since $x\in \E{P}$, it follows that $(x\cdot\varphi)(i)\in P_2$
for some $i\in\ZZ^+$. If $i>0$, then  $\varphi(i-1)=(x\cdot\varphi)(i)\in P_2$
which implies $y\in \E{P}$. If $i=0$, then $x\in P_2$ and by \eqref{eq:ip1}
applied to index pair $P$ we get $y\in P_2\subset \E{P}$.
This proves \eqref{eq:ip1} for $P^{**}$.
Consider in turn $x\in S\cup \E{P}$ and assume there exists a $y\in \Pi_{\cV}(x)\setminus(S\cup \E{P})$.
If $x\in S$ then the $\cV$-compatibility of $S$ excludes $y=x^\star$.
Therefore, we have $y\in\cl x\subset\cl S$ and by \eqref{eq:cls-s}
$y\in \E{P}$, a contradiction. This shows that $x\in \E{P}$
and proves \eqref{eq:ip2} for $P^{**}$.
Now, observe that
$S\cap \E{P}=\emptyset$ or equivalently  $(S\cup \E{P})\setminus \E{P}=S$,
which proves \eqref{eq:ip3} and the saturation property for $P^{**}$.
\qed

\subsection{Semi-equal index pairs.}
We write $P\subset Q$ for index pairs $P$, $Q$ meaning $P_i\subset Q_i$ for $i=1,2$.
We say that index pairs $P$, $Q$ of $S$
are {\em semi-equal} if $P\subset Q$ and either $P_1=Q_1$ or
$P_2=Q_2$. For semi-equal pairs $P$, $Q$, we write
\[
  A(P,Q):=\begin{cases}
             Q_1\setminus P_1 & \text{ if $P_2=Q_2$,}\\
             Q_2\setminus P_2 & \text{ if $P_1=Q_1$.}
          \end{cases}
\]

\begin{lem}
\label{lem:semi-equal}
If $P\subset Q$ are semi-equal index pairs of $S$, then $A(P,Q)$ is $\cV$-compatible.
\end{lem}
\proof
   Assume first that $P_1=Q_1$. Let $V\in\cV$ and let $x\in V$. It is enough to show that
$x\in Q_2\setminus P_2$ if and only if $x^\star\in Q_2\setminus P_2$.
Let $x\in Q_2\setminus P_2$. Since $x^\star\in\Pi_{\cV}(x)$,
we get from  \eqref{eq:ip1} applied to $P$ that $x^\star\in P_1=Q_1$.
Thus, by applying  \eqref{eq:ip1} to $Q$ we get $x^\star\in Q_2$.
Since $x\in\cl x^\star$, we get $x^\star\not\in P_2$, because otherwise $x\in P_2$.
Hence, we proved that $x\in Q_2\setminus P_2$ implies
$x^\star\in Q_2\setminus P_2$.
Let now $x^\star\in Q_2\setminus P_2$.
Then  $x\in\cl x^\star\subset Q_2$. Since $x^\star\in Q_2\subset Q_1=P_1$ and $x^\star\not\in P_2$ we
get from \eqref{eq:ip1} that $x^\star\in P_2$.
Hence, we also proved that $x^\star\in Q_2\setminus P_2$ implies
$x\in Q_2\setminus P_2$.

Consider in turn the case $P_2=Q_2$.
Again, it is enough to show that
$x\in Q_1\setminus P_1$ if and only if $x^\star\in Q_1\setminus P_1$.
Let $x\in Q_1\setminus P_1$.
Then  $x\not\in P_2=Q_2$ and by \eqref{eq:ip1} we get $x^\star\in Q_1$.
Also $x^\star\not\in P_1$, because otherwise $x\in\cl x^\star\subset P_1$.
Hence, $x\in Q_1\setminus P_1$ implies $x^\star\in Q_1\setminus P_1$.
Finally, let $x^\star\in Q_1\setminus P_1$. We have $x\in \cl x^\star\subset Q_1$.
We cannot have $x\in P_1$. Indeed, since $x^\star\not\in P_1$,  assumption $x\in P_1$ implies $x\in P_2=Q_2$.
But then $x^\star\in Q_2=P_2\subset P_1$, a contradiction. Thus, $x^\star\in Q_1\setminus P_1$
implies $x\in Q_1\setminus P_1$.
\qed

\begin{lem}
\label{lem:invA-empty-zeroSpace}
Assume $A\subset X$ is proper, $\cV$-compatible and $\Inv A=\emptyset$. Then  $A$ is a zero space.
\end{lem}
\proof
Consider a cycle
$
   x_n \adhl_{\cV} x_{n-1} \adhl_{\cV} \cdots \adhl_{\cV} x_0=x_n
$
in $A$ for some $n\geq 1$.
Then
$
  \varphi:\ZZ\ni i\mapsto x_{i \bmod{n}}\in A
$
is a solution in $A$ which shows that $\Inv A\neq\emptyset$, a contradiction.
Thus, $\cV':=\cV_{|A}$ is acyclic. Similarly, if $V\subset A$ is a critical multivector,
then $\psi:\ZZ\ni i\mapsto V^\star\in A$ is a solution in $A$,
again contradicting $\Inv A=\emptyset$. Thus, every multivector in $\cV'$ is regular.
The thesis follows now from Theorem~\ref{thm:acyclic-gdvf}
\qed

\begin{lem}
\label{lem:semi-equal-iso}
If $P\subset Q$ are semi-equal index pairs of $S$, then $H^\kappa(P_1,P_2)$ and $H^\kappa(Q_1,Q_2)$ are isomorphic.
\end{lem}
\proof
  Assume $P\subset Q$. If $P_2=Q_2$, then $A(P,Q)=Q_1\setminus P_1$. Since $Q_1$ is closed, it is proper.
Hence, since $A(P,Q)$ is open in $Q_1$, it is also proper.
Similarly we prove that $A(P,Q)$ is proper
if $P_1=Q_1$.

Thus, $A(P,Q)$ is a Lefschetz complex and it follows from Lemma~\ref{lem:semi-equal} that
the restriction $\cV':=\cV_{|A(P,Q)}$
is a multivector field on $A(P,Q)$.

Observe that
\begin{eqnarray*}
  S\cap (Q_1\setminus P_1)&\subset& P_1\cap (X\setminus P_1)=\emptyset,\\
  S\cap (Q_2\setminus P_2)&\subset& (X\setminus Q_2)\cap Q_2=\emptyset.
\end{eqnarray*}
Hence, $S\cap A(P,Q)=\emptyset$ and $\Inv A(P,Q)\subset X\setminus S$. We also have
\begin{eqnarray*}
  P_1=Q_1&\implies& \Inv(Q_2\setminus P_2)\subset \Inv(P_1\setminus P_2)=S,\\
  P_2=Q_2&\implies& \Inv(Q_1\setminus P_1)\subset \Inv(Q_1\setminus Q_2)=S.
\end{eqnarray*}
Therefore,
$
  \Inv A(P,Q)\subset S\cap (X\setminus S)=\emptyset.
$
Thus, $\cV'$ is acyclic and from Lemma~\ref{lem:invA-empty-zeroSpace} we conclude that $A(P,Q)$ is
a zero space.

It is an elementary computation to check the following two observations.
If $P_2=Q_2$, then
$P_1\setminus P_2\subset Q_1\setminus Q_2$,
$
    Q_1\setminus Q_2=A(P,Q)\cup P_1\setminus P_2
$
and $P_1\setminus P_2=P_1\cap(Q_1\setminus Q_2)$ is closed in $Q_1\setminus Q_2$.
If $P_1=Q_1$, then $Q_1\setminus Q_2\subset P_1\setminus P_2$,
$
    P_1\setminus P_2=A(P,Q)\cup Q_1\setminus Q_2
$
and $A(P,Q)=Q_2\setminus P_2=Q_2\cap(P_1\setminus P_2)$ is closed in $P_1\setminus P_2$.
Hence, by Theorem~\ref{thm:exact-sequence-L}
applied to the pair $(Q_1\setminus Q_2,P_1\setminus P_2)$
in the case $P_2=Q_2$ and
$(P_1\setminus P_2,Q_2\setminus P_2)$ in the case $P_1=Q_1$
and the fact that $A(P,Q)$ is a zero space we conclude that $H^\kappa(Q_1\setminus Q_2)$
and $H^\kappa(P_1\setminus P_2)$ are isomorphic.
\qed

\subsection{Conley index.}
The following theorem allows us to define the {\em homology Conley index} of an isolated invariant set
$S$ as $H^\kappa(P_1,P_2)=H^\kappa(P_1\setminus P_2)$ for any index pair $P$ of $S$.
We denote it $\Con(S)$.

\begin{thm}
\label{thm:conley}
Given an isolated invariant set $S$, the pair $(\cl S,\mouth S)$ is a saturated index pair for $ S$.
If $ P$ and $ Q$ are index pairs for $ S$, then $H^\kappa(P_1,P_2)$ and $H^\kappa(Q_1,Q_2)$
are isomorphic.
\end{thm}
\proof
Obviously $\cl S$ is closed and $\mouth S$ is closed by Theorem~\ref{thm:iso-inv-set}. Thus,
to show that $(\cl S,\mouth S)$ is an index pair we need to prove properties
\eqref{eq:ip1}-\eqref{eq:ip3} of the definition of index pair.
Let $x\in\mouth S$ and $y\in\Pi_{\cV}(x)\cap\cl S$. Assume $y\not\in\mouth S$.
Then $y\in S$ and
the $\cV$-compatibility of $S$ implies that $\class{x}\neq\class{y}$.
It follows that $y\in\cl x\subset\cl \mouth S=\mouth S$,
a contradiction, which shows that $y\in\mouth S$ and proves \eqref{eq:ip1}.

Consider in turn $x\in\cl S$ such that there exists a $y\in\Pi_{\cV}(x)\setminus \cl S$.
Obviously, $x\neq y$. It must be $\class{x}=\class{y}$, because otherwise $y\in\cl x\subset\cl S$.
Since $y\not\in S$, the $\cV$-compatibility of $S$ implies that
$\class{y}\cap S=\emptyset$.
But, $x\in\class{x}=\class{y}$, therefore $x\not\in S$ and consequently $x\in\mouth S$, which proves \eqref{eq:ip2}.
Obviously, $S=\cl S\setminus\mouth S$, therefore \eqref{eq:ip3} and the saturation property also are satisfied.
This completes the proof of the first part of the theorem.

To prove the remaining assertion first observe that it is obviously satisfied if both pairs are saturated.
Thus, it is sufficient to prove the assertion in the case $Q=P^{**}$, because by Lemma~\ref{lem:saturated}
the pair $P^{**}$ is saturated. We obviously have $P^*\subset P$ and $P^*\subset P^{**}$ and each inclusion
is a semi-equality. Therefore, both inclusions induce an isomorphism in homology by
Lemma~\ref{lem:semi-equal-iso}.
\qed

We call the index pair $(\cl S,\mouth S)$ {\em canonical},
because it minimizes both $P_1$ and $P_1\setminus P_2$.
Indeed, for any index pair $\cl S\subset P_1$, because $P_1$ is closed and by \eqref{eq:ip3}
$S\subset P_1\setminus P_2$. In the case of the canonical index pair both inclusions are equalities.
Note that in the  case of the classical Conley index there is no natural choice of a canonical index pair.
Since  $S=\cl S\setminus\mouth S$, we see that  $\Con(S)\cong H^\kappa(S)$.
Thus, in our combinatorial setting Theorem~\ref{thm:conley} is actually not needed to define
the Conley index. However, the importance of Theorem~\ref{thm:conley} will become clear in
the proof of Morse inequalities via the following corollary.

\begin{cor}
\label{cor:conley}
If $(P_1,P_2)$ is an index pair of an isolated invariant set $S$, then
\begin{equation}
\label{eq:conley}
   p_S(t)+p_{P_2}(t)=p_{P_1}(t)+(1+t)q(t),
\end{equation}
where $q(t)$ is a polynomial with non-negative coefficients.
Moreover, if $H^\kappa(P_1)=H^\kappa(P_2)\oplus H^\kappa(S)$,
then $q(t)=0$.
\end{cor}
\proof
   By applying Theorem~\ref{thm:exact-sequence-L} to the pair $P_2\subset P_1$ of closed subsets of $X$
and Theorem~\ref{thm:exact-sequence}
to the resulting exact sequence we obtain
\[
   p_{P_1\setminus P_2}(t)+p_{P_2}(t)=p_{P_1}(t)+(1+t)q(t)
\]
for some polynomial $q$ with non-negative coefficients.
The conclusion follows now from Theorem~\ref{thm:conley}
by observing that $H^\kappa(P_1\setminus P_2)=H^\kappa(S)$.
\qed

In the case of an isolated invariant set $S$ we call the polynomial
$p_S(t)$ the {\em Conley polynomial} of $S$.
We also define the $i$th Conley coefficient of $S$ as
\[
   c_i(S):=\rank H^\kappa_i(S).
\]

\subsection{Additivity of the Conley index.}
Let $S$ be an isolated invariant set and assume $S_1,S_2\subset S$ are also isolated invariant sets.
We say that $S$ {\em decomposes} into  $S_1$ and $S_2$ if $S_1\cap S_2=\emptyset$ and
every full solution $\varrho$ in $S$ is either a solution in $S_1$ or in $S_2$.
\begin{thm}
\label{thm:additivity}
Assume an isolated invariant set $S$ decomposes into the union of two its isolated invariant subsets $S_1$ and $S_2$.
Then  $\Con(S)=\Con(S_1)\oplus\Con(S_2)$.
\end{thm}
\proof
First observe that the assumptions imply that $S=S_1\cup S_2$.
We will show that $S_1$ is closed in $S$.
Let $x\in \cl_S S_1$ and assume $x\not\in S_1$. Then  $x\in S_2$.
Choose $y\in S_1$ such that $x\in\cl_S y$.
The $\cV$-compatibility of $S_1$ implies that
$\class{x}\neq\class{y}$, hence $x\in\Pi_{\cV}(y)$.
Select $\gamma\in\Sol(y,S_1)$ and $\varrho\in\Sol(x,S_2)$.
Then $\varphi:=\gamma^-\cdot \varrho^+$ is a solution in $S$
but neither in $S_1$ nor in $S_2$, which contradicts the assumption that
$S$ decomposes into  $S_1$ and $S_2$. Hence, $S_1$ is closed.
Similarly we prove that $S_2$ is closed. The conclusion follows now from
Proposition~\ref{prop:hom-union}.
\qed

\section{Attractors and repellers.}
\label{sec:attractors-repellers}

In this section we  define attractors and repellers and study attractor-repeller pairs
needed to prove Morse equation and Morse inequalities.
In the whole section we assume that the Lefschetz complex $X$ is invariant with respect
to a given combinatorial multivector field $\cV$ on $X$.

\subsection{Attractors.}
We say that a $\cV$-compatible $N\subset X$ is a {\em trapping region}
if $\Pi_{\cV}(N)\subset N$. This is easily seen to be equivalent
to the requirement that $N$ is $\cV$-compatible and for any $x\in N$
\begin{equation}
\label{eq:trapping}
\Sol^+(x,X) = \Sol^+(x,N).
\end{equation}
We say that $ A$ is an {\em attractor}
if there exists a trapping region $N$ such that
$ A=\Inv N$.

\begin{thm}
\label{thm:attractor}
The following conditions are equivalent:
\begin{itemize}
  \item[(i)] $ A$ is an attractor,
  \item[(ii)] $ A$ is invariant and closed,
  \item[(iii)] $ A$ is isolated invariant and closed,
  \item[(iv)] $ A$ is isolated invariant, closed and a trapping region.
\end{itemize}
\end{thm}
\proof
  Assume (i). Let $N$ be a trapping region such that $A=\Inv N$.
By \eqref{eq:inv2} we have $\Inv A=\Inv \Inv N=\Inv N=A$, hence $A$ is
invariant. To see that $A$ is closed take $x\in\cl A$ and $y\in A$ such
that $x\in\cl y\subset\cl y^\star$.
Note that $y^\star\in A$, because by Proposition~\ref{prop:inv-set}
the set $A$, as invariant, is $\cV$-compatible.
If $[x]=[y]$, then for the same reason $x\in A$.
Hence, assume $[x]\neq [y]=[y^\star]$.
Then  $x\in\Pi_{\cV}(y^\star)\subset\Pi_{\cV}(N)\subset N$.
From the $\cV$-compatibility of $N$ we get that also $x^\star\in N$.
Let $\varphi\in\Sol^+(x^\star)$.
Since $N$ is a trapping region, we see
that $\varphi\in\Sol^+(x^\star,N)$.
It follows that $x\in\Inv^+N$.
Since $y^\star\in A=\Inv N\subset\Inv^-N$, we can take $\psi\in\Sol^-(y^\star,N)$.
Then   $\psi\cdot \nu(x)\in\Sol^-(x^\star,N)$. This shows that $x\in\Inv^-N$.
By \eqref{eq:inv1} we have $x\in\Inv N=A$. This proves that $A$ is closed and shows
that (i) implies (ii).
If (ii) is satisfied, then $A$, as closed, is proper.
Hence, we get from Theorem~\ref{thm:iso-inv-set} that (ii) implies (iii).
Assume in turn (iii).
Let $x\in A$ and let $y\in\Pi_{\cV}(x)$.
If $[y]=[x]$, then
$y\in A$ by the $\cV$-compatibility of $A$. If $[y]\neq [x]$,
then $x=x^\star$,  $y\in\cl x^\star\subset \cl A=A$. Hence $\Pi_{\cV}(A)\subset A$,
which proves (iv). Finally, observe that (i) follows immediately from (iv).
\qed

\subsection{Repellers.}
We say that a $\cV$-compatible $N\subset X$ is a {\em backward trapping region}
if $\Pi_{\cV}^{-1}(N)\subset N$. This is easily seen to be equivalent
to the requirement that $N$ is $\cV$-compatible and  for any $x\in N$
\begin{equation}
\label{eq:b-trapping}
\Sol^-(x,X) = \Sol^-(x,N).
\end{equation}
We say that $R$ is a {\em repeller}
if there exists a backward trapping region $N$ such that $ R=\Inv N$.
The following proposition is straightforward.
\begin{prop}
\label{prop:f-b-trapping}
A subset $N\subset X$ is a trapping region if and only if $X\setminus N$
is a backward trapping region.
\qed
\end{prop}

\begin{thm}
\label{thm:repeller}
The following conditions are equivalent:
\begin{itemize}
   \item[(i)] $ R$ is a repeller,
   \item[(ii)] $ R$ is invariant and open,
   \item[(iii)] $ R$ is isolated invariant and open,
   \item[(iv)] $ R$ is isolated invariant, open and a backward trapping region.
\end{itemize}
\end{thm}
\proof
  Assume (i). Let $N$ be a backward trapping region such that $R=\Inv N$.
By \eqref{eq:inv2} we have $\Inv R=\Inv \Inv N=\Inv N=R$, hence $R$ is
invariant. To see that $R$ is open
we will prove that $X\setminus R$ is closed.
For this end take $x\in\cl(X\setminus R)$ and $y\in X\setminus R$ such that $x\in\cl y\subset\cl y^\star$.
Note that $y^\star\in X\setminus R$, because
the set $X\setminus R$ is $\cV$-compatible as a complement of
an invariant set which, by Proposition~\ref{prop:inv-set}, is $\cV$-compatible.
If $[x]=[y]$, then for the same reason $x\in X\setminus R$.
Hence, assume $[x]\neq [y]=[y^\star]$.
Then  $x\in\Pi_{\cV}(y^\star)$.
Let $\psi\in\Sol^+(x^\star)$ and $\varphi\in\Sol^-(y^\star)$.
Then  $\gamma:=\varphi\cdot \nu^-(x)\cdot \psi\in\Sol(y^\star)$.
Since $y^\star\not\in R$, there exists an $i\in\ZZ$ such that $\gamma(i)\not\in N$.
By Proposition~\ref{prop:f-b-trapping} we see that
$X\setminus N$ is a trapping region. Hence, $\gamma(j)\not\in N$ for all $j\geq i$.
In particular, $\psi(j)\not\in N$ for large $j$.
Since $\psi\in\Sol^+(x^\star)$ is arbitrary, it follows that $x\not\in\Inv^+N$ and consequently
$x\not\in\Inv N=R$.
Hence, $X\setminus R$ is closed, that is $R$ is open.
This proves (ii).
If (ii) is satisfied, then $R$, as open, is proper.
Hence, we get from Theorem~\ref{thm:iso-inv-set} that (ii) implies (iii).
Assume (iii).
Let $x\in R$ and let $y\in\Pi_{\cV}^{-1}(x)$. If $[y]=[x]$, then
$y\in R$ by the $\cV$-compatibility of $R$. If $[y]\neq [x]$,
then $y=y^\star$ and  $x\in\cl y^\star$.
It follows from Proposition~\ref{prop:f-topology} that
$y=y^\star\in\opn x\subset R$, because $x\in R$ and $R$ is open.
Hence, $\Pi_{\cV}^{-1}(R)\subset R$,
which proves (iv). Finally, observe that (i) follows immediately from (iv).
\qed

\subsection{Recurrence and basic sets.}
Let $A\subset X$ be $\cV$-compatible.
We write $x\rpthVA{A} y$ if there exists a path of $\cV$ in $A$ from $x^\star$ to $y^\star$ of length at least two.
We write $x\epthVA{A} y$ if $x\rpthVA{A} y$ and $y\rpthVA{A} x$.
We drop $\cV$ in this notation if $\cV$ is clear from the context.
Also, we drop $A$ if $A=X$.
We say that $A$ is {\em weakly recurrent} if for every $x\in A$ we have $x\epthA{A} x$.
It is {\em strongly recurrent} if for any $x,y\in A$
we have $x\rpthA{A} y$, or equivalently for any $x,y\in A$
there is $x\epthA{A} y$.
Obviously, every strongly recurrent set is also weakly recurrent.

The {\em chain recurrent set of $ X$}, denoted $\CR(X)$, is the maximal weakly recurrent subset of $X$.
It is straightforward to verify that
\[
  \CR(X):=\setof{x\in X\mid x\epthV x}.
\]
Obviously, the relation $\epthV$ restricted to
$\CR(X)$ is an equivalence relation.
By a {\em basic set} of $\cV$ we mean an equivalence class of $\epthV$ restricted to
$\CR(X)$.
\begin{thm}
\label{thm:basic set}
Every basic set is a strongly recurrent isolated invariant set.
\end{thm}
\proof
Let $B$ be a basic set. Obviously, it is $\cV$-compatible and invariant.
By its very definition it is also strongly recurrent.
Thus, by Theorem~\ref{thm:iso-inv-set} we only need to prove that $B$ is proper.
For this end  we will verify condition \eqref{eq:proper-sets} of Proposition~\ref{prop:proper-sets}.
Take $x,z\in B$, $y\in X$ and assume $x\in\cl y$ and $y\in \cl z$. Then also  $x\in\cl y^\star$ and $y\in \cl z^\star$.
If $\class{x}=\class{y^\star}$ or $\class{y}=\class{z^\star}$, then $y\in B$ by the $\cV$-compatibility of $B$.
Thus, assume $\class{x}\neq\class{y^\star}$ and $\class{y}\neq\class{z^\star}$.
It follows that $x\in\Pi_{\cV}(y^\star)$ and  $y\in\Pi_{\cV}(z^\star)$ and consequently $y^\star\rpthA{B} x^\star$
and $z^\star\rpthA{B} y^\star$. By the definition of the basic set also $x^\star\rpthA{B} z^\star$.
It follows that the basic sets of $x$, $y$ and $z$ coincide and consequently $y\in B$.
\qed

\subsection{Limit sets.}

Let $\varrho:\ZZ\to X$ be a full solution.
Recall that $[A]_\cV^+$ denotes the intersection of $\cV$-compatible sets containing $A$.
We define the {\em $\alpha$ and $\omega$ limit sets of $\varrho$}
by
\begin{eqnarray*}
  \alpha(\varrho)&:=&\bigcap_{k\geq 0}\Inv [\im \sigma_-^k\varrho^-]^+_{\cV},\\
  \omega(\varrho)&:=&\bigcap_{k\geq 0}\Inv [\im \sigma_+^k\varrho^+]^+_{\cV}.
\end{eqnarray*}

\begin{prop}
\label{prop:varrho-inv}
The sets $\alpha(\varrho)$ and $\omega(\varrho)$ are
invariant, $\cV$-compatible and strongly recurrent.
Moreover, $\alpha(\varrho)\cap\im \varrho\neq\emptyset\neq\omega(\varrho)\cap\im \varrho$.
\end{prop}
\proof
Let $\varphi_k:=\sigma_+^k\varrho^+$ and $A_k:=\Inv[\im\varphi_k]^+_{\cV}$.
Obviously, the sets $A_k$ are invariant and $\cV$-compatible.
It is straightforward to observe that $A_{k+1}\subset A_{k}$.
Hence, since $X$ is finite, there exists a $p\in\ZZ^+$ such that
$A_{p}=\bigcap_{k\geq 0} A_k=\omega(\varrho)$. This proves that
$\omega(\varrho)$ is invariant and $\cV$-compatible.
To see that $\omega(\varrho)$ is strongly recurrent, take $x,y\in\omega(\varrho)=A_p$.
Then  $x^\star,y^\star\in \im \varphi_p\cap A_p$. Let $m,n\geq p$ be such that $x^\star=\varphi(m)$
and $y^\star=\varphi(n)$. If $n<m$, then $\varphi_{|[n,m]}$ is a path in $A_p$ from $x^\star$ to $y^\star$.
Thus, assume $m\leq n$
and take $q>m$. Since $A_q=A_p$, we see that $y^\star\in A_q$. Thus, there exists a $k\geq q$
such that $y^\star=\varphi(k)$. Hence, $\varphi_{|[m,k]}$ is a path in $A_p$ from $x^\star$ to $y^\star$.
Therefore, $A_p=\omega(\varrho)$ is strongly recurrent.

To see that $\alpha(\varrho)\cap\im \varrho\neq\emptyset$
observe that there exist $m,n\in\ZZ^+$, $m<n$,
such that $\varphi_p(m)=\varphi_p(n)$.
By Proposition~\ref{prop:solution-star} we may assume without loss of generality that $\varphi_p(m)=\varphi_p(m)^\star$.
Consider
\[
   \psi:\ZZ\ni i \mapsto \varphi_p(m+i\bmod{m-n})\in\im\varphi_p.
\]
We have $\psi\in\Sol(\varphi_p(m)^\star,[\im\varphi_p]^+_\cV)$,
therefore $\varphi_p(m)\in \Inv [\im\varphi_p]^+_{\cV}=A_p=\omega(\varrho)$.
Obviously, $\varphi_p(m)\in\im\varphi$, which
proves that
$\omega(\varrho)\cap\im\varrho\neq\emptyset$.
The proof concerning $\alpha(\varrho)$ is analogous.
\qed

\subsection{Attractor-repeller pairs}
Since by Theorem~\ref{thm:attractor} an attractor is in particular a trapping region,
it follows from Proposition~\ref{prop:f-b-trapping} that if $A$
is an attractor, then $\Inv(X\setminus A)$ is a repeller. It is called
the {\em dual repeller of $A$} and is denoted by $A^\star$.
Similarly, if $R$ is a repeller, then $\Inv(X\setminus R)$ is an attractor,
called the {\em dual attractor of $R$} and denoted by $R^\star$.

For $A,A'\subset X$ define
\[
  C(A',A):=\setof{ x\in X\mid\exists_{\varrho\in \Sol(x^\star)}\; \alpha(\varrho)\subset A',\;\omega(\varrho)\subset A}.
\]
\begin{prop}
\label{prop:no-connection}
Assume $A$ is an attractor. Then  there is no heteroclinic connection
running from  $A$ to $A^\star$, that is $C(A,A^\star)=\emptyset$.
\qed
\end{prop}
The pair $( A, R)$ of subsets of $ X$ is said to be an {\em attractor-repeller pair in $ X$}
if $A$ is an attractor, $R$ is a repeller, $A=R^\star$ and $R=A^\star$.
\begin{thm}
\label{thm:attractor-repeller}
  The pair $(A,R)$ is an attractor-repeller pair in $ X$ if and only if the following four conditions are satisfied:
  \begin{itemize}
     \item[(i)] $A$ is an attractor and $R$ is a repeller,
     \item[(ii)] $ A\cap R=\emptyset$,
     \item[(iii)] for every $x\not\in A$ and $\varrho\in\Sol(x)$ we have $\alpha(\varrho)\subset R$,
     \item[(iv)] for every $x\not\in R$ and $\varrho\in\Sol(x)$ we have $\omega(\varrho)\subset A$.
  \end{itemize}
\end{thm}
\proof
First assume that $(A,R)$ is an attractor-repeller pair.
Then, obviously, (i) is satisfied. Since
$R=A^\star=\Inv(X\setminus A)\subset X\setminus A$, we see that (ii) holds.
In order to prove (iii) take  $x\not\in A$
and $\varrho\in\Sol(x)$. Since $X\setminus A$ is a backward
trapping region, we see that $\im\varrho^-\subset X\setminus A$.
Moreover, $[\im\varrho^-]^+_{\cV}\subset X\setminus A$, because $X\setminus A$,
as a backwards trapping region, is $\cV$-compatible.
Hence, $\alpha(\varrho)\subset\Inv [\im\varrho^-]^+_{\cV}\subset \Inv(X\setminus A)=A^\star=R$.
Similarly we prove (iv).

Assume in turn that (i)-(iv) hold.
We will show that $R=A^\star$. Indeed, by (ii) $R\subset X\setminus A$,
hence $R=\Inv R\subset \Inv(X\setminus A)=A^\star$. To show
the opposite inclusion take $x\in A^\star$ and assume $x\not\in R$.
Let $\varrho\in\Sol(x^\star,X\setminus A)$. By (iv) we have $\omega(\varrho)\subset A$.
By Proposition~\ref{prop:varrho-inv} we can select a $y\in\omega(\varrho)\cap\im\varrho$.
It follows that $y\in A\cap\im\varrho$. This contradicts $\varrho\in\Sol(x^\star,X\setminus A)$ and proves that $R=A^\star$.
Similarly we prove that $A=R^\star$.
\qed

The following corollary is an immediate consequence of Theorem~\ref{thm:attractor-repeller}.
\begin{cor}
\label{cor:attractor-repeller}
Assume $A$ is an attractor. Then $(A,A^\star)$ is an attractor-repeller pair.
In particular, $A^{**}=A$.
\qed
\end{cor}

\section{Morse decompositions, Morse equation and Morse inequalities}
\label{sec:morse}

In this section we define Morse decompositions and we prove the Morse equation
and Morse inequalities for Morse decompositions.

\subsection{Morse decompositions}
Let $\PP$ be a finite set.
The collection $\cM=\setof{M_r\mid r\in \PP}$ is called a {\em Morse decomposition of $X$} if
there exists a partial order $\leq$ on $\PP$ such that
the following three conditions are satisfied:
\begin{itemize}
   \item[(i)] $\cM$ is a family of mutually disjoint, isolated invariant subsets
              of $X$,
   \item[(ii)] for every full solution $\varphi$ in $X$ there exist
               $r,r'\in \PP$, $r\leq r'$,  such that $\alpha(\varphi)\subset M_{r'}$,
               $\omega(\varphi)\subset M_{r}$,
   \item[(iii)] if for a full solution $\varphi$ and $r\in \PP$ we have
                $\alpha(\varphi)\cup\omega(\varphi)\subset M_{r}$, then $\im\varphi\subset M_r$.
\end{itemize}
Note that since the sets in $\cM$ are mutually disjoint, the indices $r,r'$ in (ii) are
determined uniquely. A partial order on $\PP$ which makes $\cM=\setof{M_r\mid r\in \PP}$
a Morse decomposition of $X$ is called an {\em admissible} order.
Observe that the intersection of admissible orders is an admissible order and any extension
of an admissible order is an admissible order. In particular, every Morse decomposition
has a unique minimal admissible order as well as an admissible linear order.

The following proposition is straightforward.
\begin{prop}
\label{prop:C-p-prim-p}
Let  $\cM=\setof{M_r\mid r\in \PP}$ be a Morse decomposition of $X$ and let $r,r'\in \PP$. Then
\begin{itemize}
   \item[(i)] $C(M_{r'},M_r)$ is $\cV$-compatible,
   \item[(ii)] $C(M_r,M_r)=M_r$,
   \item[(iii)] $r'<r$ implies $C(M_{r'},M_r)=\emptyset$.
\end{itemize}
\qed
\end{prop}

Let $\cB$ be the family of all basic sets of $\cM$.
For $B_1,B_2\in\cB$ we write $B_1\leqV B_2$ if there exists a solution
$\varrho$ such that $\alpha(\varrho)\subset B_2$ and $\omega(\varrho)\subset B_1$.
It follows easily from the definition of the basic set that $\leqV$ is a partial order on $\cB$.

\begin{thm}
\label{thm:basic-sets}
The family $\cB$ of all basic sets of $\cV$ is a Morse decomposition of $X$.
\end{thm}
\proof
Obviously, two different basic sets are always disjoint.
By Theorem~\ref{thm:basic set} each element of $\cB$ is an isolated invariant set.
Thus, condition (i) of the definition of Morse decomposition is satisfied.
To prove (ii) consider a full solution $\varrho$. By Proposition~\ref{prop:varrho-inv} both
$\alpha(\varrho)$ and $\omega(\varrho)$ are strongly recurrent. Thus, each is contained
in a basic set, which proves (ii).
Assume in turn that both $\alpha(\varrho)$ and $\omega(\varrho)$ are contained in the same basic set $B$.
Fix $y\in\im\varrho$ and choose $x\in\alpha(\varrho)$ and $z\in\omega(\varrho)$.
Then $x\rpthV y$ and $y\rpthV z$. But also $z\rpthV x$, because $x,z\in B$ and $B$ is strongly recurrent.
Thus, $x\epthV y$ and consequently $y\in B$. This shows that $\im\varrho\subset B$ and proves (iii).
Finally, observe that $\leqV$ is obviously and admissible partial order on $\cB$.
Therefore, $\cB$ is a Morse decomposition of $X$.
\qed

Given two Morse decompositions $\cM$, $\cM'$ of $X$, we write $\cM'\sqsubset\cM$
if for each $M'\in\cM'$ there exists an $M\in\cM$ such that $M'\subset M$.
Then, we say that $\cM'$ is {\em finer} than $\cM$. The relation $\sqsubset$
is easily seen to be a partial order on the collection of all Morse decompositions of $X$.

\begin{thm}
\label{thm:finest}
For any Morse decomposition $\cM=\{M_r\mid r\in\PP\}$ of $X$ we have $\cB\sqsubset\cM$.
Thus, the family of basic sets is the unique, finest Morse decomposition of $X$.
\end{thm}
\proof
  Let $B$ be a basic set and let $x\in B$. Since $B$ is invariant, we may choose a $\varrho\in\Sol(x^\star,B)$
and $r\in\PP$ such that $\alpha(\varrho)\subset M_r$. Since $\alpha(\varrho)\subset B$, we may choose
a $y=y^\star\in \im\varrho\cap B\cap M_r$. We will show that $B\subset M_r$. For this end take a $z\in B$.
Since $B$ is strongly recurrent, we may construct a path $\varphi$ from $y$ to $y$ through $z$.
Then $\gamma:=\sigma_-\varrho^-\cdot\varphi\cdot\sigma_+\varrho^+$ is a solution through $z$ and obviously
$\alpha(\gamma)=\alpha(\varrho)\subset M_r$ and $\omega(\gamma)=\omega(\varrho)\subset M_r$.
It follows that $z\in \im \gamma\subset M_r$. This proves our claim.
\qed

\subsection{Morse sets.}

Given $I\subset \PP$ define the {\em Morse set }of $I$ by
\begin{equation}
\label{eq:M-I}
   M(I):=\bigcup_{r,r'\in I} C(M_{r'},M_r).
\end{equation}

\begin{thm}
\label{thm:M-I}
  The set $M(I)$ is an isolated invariant set.
\end{thm}
\proof
To prove that $M(I)$ is invariant take $x\in M(I)$ and let $r,r'\in \PP$
be such that $x\in C(M_{r'},M_r)$.
Choose $\varrho\in\Sol(x^\star)$ such that $\alpha(\varrho)\subset M_{r'}$
and $\omega(\varrho)\subset M_{r}$. Let $k\in\ZZ$ and let $y:=\varrho(k)$.
Since either $y^\star=\varrho(k)$ or $y^\star=\varrho(k+1)$, we see that
$\rho\in\Sol(y^\star)$. It follows that $y\in C(M_{r'},M_r)\subset M(I)$. Therefore,
$\rho\in\Sol(x^\star,M(I))$ and $x\in\Inv M(I)$.

By Theorem~\ref{thm:iso-inv-set} it suffices
to prove that $M(I)$ is proper.
For this end  we will verify condition \eqref{eq:proper-sets}
in Proposition~\ref{prop:proper-sets}.
Take $x,z\in M(I)$ and assume $y\in X$
is such that $x\in\cl y$ and $y\in\cl z$. If $[x]=[y]$ or $[y]=[z]$, then $y\in M(I)$,
because by Proposition~\ref{prop:inv-set} the set  $M(I)$ as invariant is $\cV$-compatible.
Hence, consider the case
$[x]\neq [y]$ and $[y]\neq [z]$. It follows that $x\in\Pi_{\cV}(y^\star)$
and $y\in\Pi_{\cV}(z^\star)$. Let $\gamma\in\Sol(x^\star)$ be such that
$\omega(\gamma)\subset M_{r}$
for some $r\in I$.
Similarly, let $\varrho\in\Sol(z^\star)$ be such that
$\alpha(\varrho)\subset M_{s}$ for some $s\in I$.
Let $\chi:=\varrho^-\cdot\nu(y)\cdot\nu^-(x)\cdot\gamma^+$.
We easily verify that $\alpha(\chi)=\alpha(\varrho)\subset M_s$
and $\omega(\chi)=\omega(\gamma)\subset M_r$. This shows that $y\in M(I)$.
Thus, $M(I)$ is proper.
\qed

\begin{thm}
\label{thm:M-A}
If $I$ is a lower set in $\PP$, then $M(I)$ is an attractor in $ X$.
\end{thm}
\proof
   In view of Theorem~\ref{thm:attractor} and Theorem~\ref{thm:M-I}
we only need to prove that $M(I)$ is closed. Hence, take $x\in\cl M(I)$.
We need to show that $x\in M(I)$.
Let $y\in M(I)$ be such that $x\in\cl y$. If $[x]=[y]$, then $x\in M(I)$,
because $M(I)$ as invariant is $\cV$ compatible. Thus, assume $[x]\neq [y]$. Then  $x\in\Pi_{\cV}(y^\star)$.
Let  $\gamma\in\Sol(y^\star)$ be such that $\alpha(\gamma)\subset M_{r}$ for some
$r\in I$. Choose also a $\varrho\in\Sol(x^\star)$ and let $s\in \PP$ be such that
$\omega(\varrho)\subset M_s$. Then  $\chi:=\gamma^-\cdot\nu^-(x)\cdot\varrho^+\in\Sol(x^\star)$
 and $\alpha(\chi)=\alpha(\gamma)\subset M_r$,
$\omega(\chi)=\omega(\varrho)\subset M_s$. It follows from the definition
of Morse decomposition that $s\leq r$. Since $r\in I$ and $I$ is a lower set,
we get $s\in I$. This implies that $x\in M(I)$.
\qed

We also have a dual statement. We skip the proof, because
it is similar to the proof of Theorem~\ref{thm:M-A}.
\begin{thm}
\label{thm:M-R}
If $I$ is an upper set in $\PP$, then $M(I)$ is a repeller in $X$.
\end{thm}

\begin{prop}
\label{prop:forw-back-trap}
   If $N$ is a trapping region and $N'$ is a backward trapping region, then
$
   \Inv(N\cap N')=\Inv N\cap \Inv N'.
$
\end{prop}
\proof
   Obviously, the left-hand-side is contained in the right-hand-side.
To prove the opposite inclusion take $x\in \Inv N\cap \Inv N'$.
By the $\cV$-compatibi\-li\-ty of invariant sets also $x^\star\in\Inv N\cap\Inv N'$.
Let $\varrho\in\Sol(x^\star,N')$. Since $N$ is
a trapping region we have $\varrho^+\in\Sol^+(x^\star,N)$.
Hence, $x\in\Inv^+(N\cap N')$.
Similarly we prove that  $x\in\Inv^-(N\cap N')$. It follows that $x\in\Inv(N\cap N')$.
\qed

For a lower set $I$ in $\PP$ set
\[
  N(I):=X\setminus M(I)^\star.
\]
Observe that by Corollary~\ref{cor:attractor-repeller}  we have
\begin{equation}
\label{eq:InvNIisMI}
   \Inv N(I)=M(I).
\end{equation}

\begin{thm}
\label{thm:index-pair}
  If $I\subset \PP$ is convex, then $(N(I^{\leq}),N(I^{<}))$ is an index pair for $ M(I)$.
\end{thm}
\proof
  The sets $M(I^{\leq})^\star$, $M(I^{<})^\star$ are repellers.
Hence, by Theorem~\ref{thm:repeller} they are open.
It follows that $N(I^{\leq})$, $N(I^{<})$ are closed.
To verify property \eqref{eq:ip1} take $x\in N(I^{<})$, $y\in \Pi_{\cV}(x)\cap N(I^{\leq})$ and
assume that $y\not\in N(I^{<})$. Then $y\in M(I^{<})^\star$.
It follows that
$
x\in\Pi_{\cV}^{-1}(M(I^{<})^\star)\subset M(I^{<})^\star,
$
because by Theorem~\ref{thm:repeller} the set $M(I^{<})^\star$, as a repeller, is a backwards trapping region.
Thus, $x\not\in N(I^<)$, a contradiction which proves \eqref{eq:ip1}.
To prove \eqref{eq:ip2} take $x\in N(I^{\leq})$ such that $\Pi_{\cV}(x)\setminus N(I^{\leq})\neq\emptyset$
and assume that $x\not\in N(I^{<})$. Then  $x\in  M(I^{<})^\star$. Let $y\in \Pi_{\cV}(x)\setminus N(I^{\leq})$.
Then  $y\in \Pi_{\cV}(x)\cap  M(I^{\leq})^\star$. Since $M(I^{\leq})^\star$ is a repeller, it follows that $x\in M(I^{\leq})^\star$.
This contradicts $x\in N(I^{\leq})$ and proves \eqref{eq:ip2}.

Finally, to prove \eqref{eq:ip3} first observe that
\begin{equation}
\label{eq:index-pair-aux}
  N(I^{\leq})\setminus N(I^{<})=M(I^{<})^\star\setminus M(I^{\leq})^\star=M(I^{<})^\star\cap N(I^{\leq}).
\end{equation}
Since $M(I^{<})^\star$ is a backward trapping region and $N(I^{\leq})$ is a forward trapping region, we have
by \eqref{eq:index-pair-aux} and Proposition~\ref{prop:forw-back-trap}
\begin{eqnarray*}
  \Inv( N(I^{\leq})\setminus N(I^{<})) &=& \Inv M(I^{<})^\star\cap \Inv N(I^{\leq})=\\
                                 &&  M(I^{<})^\star\cap  M(I^{\leq})=\\
                                 && \Inv M(I^{\leq})\cap \Inv(X\setminus M(I^<))=\\
                                 && \Inv( M(I^{\leq})\setminus M(I^{<}))=M(I).
\end{eqnarray*}
\qed

The following corollary is an immediate consequence of Theorem~\ref{thm:index-pair}
and the definition of the lower set.

\begin{cor}
\label{cor:index-pair-attractor}
  If $I$ is a lower set (an attracting interval) in $\PP$, then $I^{\leq}=I$,
  $I^{<}=\emptyset$,  $(N(I),\emptyset)$ is an index pair for $M(I)$, $H^\kappa(N(I))=H^\kappa(M(I))$
  and $p_{M(I)}(t)=p_{N(I)}(t)$.
\qed
\end{cor}

\begin{thm}
\label{thm:ar-equation}
Assume  $A\subset X$ is an attractor and $A^\star$ is the dual repeller.
Then
\begin{equation}
\label{eq:ar-equation}
    p_A(t)+p_{A^\star}(t)=p_X(t)+(1+t)q(t)
\end{equation}
for a polynomial $q(t)$ with non-negative coefficients.
Moreover, if $q\neq 0$, then $C(A^\star,A)\neq\emptyset$.
\end{thm}
\proof
Take $\PP:=\{1,2\}$, $M_2:=A^\star$, $M_1:=A$.
Then $\cM:=\{M_1,M_2\}$ is a Morse decomposition of $X$. Take $I:=\{2\}$.
By Proposition~\ref{prop:lower-integer} we have $I^\leq=\{1,2\}$, $I^<=\{1\}$.
It follows that $M(I^\leq)=X$, $M(I^<)=M(\{1\})=A$, $N(I^\leq)=X\setminus X^\star=X$.
Thus, we get from Corollary~\ref{cor:index-pair-attractor} that
\begin{equation}
\label{eq:ar-equation-2}
   p_{N(I^\leq)}(t)=p_X(t), \quad p_{N(I^<)}(t)=p_{M(I^<)}(t)=p_A(t).
\end{equation}
By Theorem~\ref{thm:index-pair} the pair $(N(I^\leq),N(I^<))$ is an index pair for $M(I)=A^\star$.
Thus, by substituting $P_1:=N(I^\leq)$, $P_2:=N(I^<)$, $S:=A^\star$ into \eqref{eq:conley}
in Corollary~\ref{cor:conley}
we get \eqref{eq:ar-equation} from \eqref{eq:ar-equation-2}.
By Proposition~\ref{prop:no-connection} we have $C(A,A^\star)=\emptyset$.
If also $C(A^\star,A)=\emptyset$, then $X$ decomposes into $A$, $A^\star$
and by Theorem~\ref{thm:additivity}
\[
   H^\kappa(P_1)=\Con(X)=\Con(A)\oplus\Con(A^\star)=H^\kappa(P_2)\oplus H^\kappa(A^\star).
\]
Thus, $q=0$ by Corollary~\ref{cor:conley}.
It follows that $q\neq 0$ implies $C(A^\star)\neq\emptyset$.
\qed
\subsection{Morse equation.}
\label{sec:morse-equation}

We begin by the observation that an isolated invariant set $S$ as a proper and $\cV$-compatible subset
of $X$ in itself may be considered a Lefschetz complex
with a combinatorial multivector field being the restriction $\cV_{|S}$.
In particular, it makes sense to consider attractors, repellers and Morse decompositions of $S$
with respect to $\cV_{|S}$ and the results of the preceding sections
apply.

\begin{thm}
\label{thm:Morse-equation}
   Assume $\PP=\setof{1,2,\ldots n}$ is ordered by the linear order  of natural numbers.
Let $\cM:=\setof{ M_p\mid p\in \PP}$ be a Morse decomposition of an isolated invariant set $S$
and let $A_i:=M(\{i\}^\leq)$ and $A_0:=\emptyset$.
Then  $(A_{i-1},M_i)$ is an attractor-repeller pair in $A_i$.
Moreover,
\begin{equation}
\label{eq:Morse-equation}
  \sum_{i=1}^np_{M_i}(t)=p_S(t)+(1+t)\sum_{i=1}^nq_i(t)
\end{equation}
for some polynomials $q_i(t)$ with non-negative coefficients
and such that $q_i(t)\neq 0$ implies $C(M_i,A_{i-1})\neq\emptyset$
for $i=2,3,\ldots n$.
\end{thm}
\proof
Fix $i\in\PP$.
Obviously, $A_{i-1}\subset A_i$. The set $A_{i-1}$, as an attractor in $S$, is closed in $S$.
Thus, it is closed in $A_i$. It follows that $A_{i-1}$ is an attractor in $A_i$.
The verification that $M_i$ is the dual repeller of $A_{i-1}$ in $A_i$
is straightforward. By applying  Theorem~\ref{thm:ar-equation} to the attractor-repeller pair
$(A_{i-1},M_i)$ in $A_i$ we get
\begin{equation}
\label{eq:Morse-equation-2}
p_{M_{i}}(t) + p_{A_{i-1}}(t) = p_{A_{i}}(t) + (1+t)q_{i}(t)
\end{equation}
for a polynomial
$q_{i}$ with  non-negative coefficients.
Since $p_{A_{0}}(t)=0$ and $A_{n}=S$,
summing \eqref{eq:Morse-equation-2} over $i\in\PP$  and  substituting  $q:=\sum_{i=1}^n q_{i},$
we get \eqref{eq:Morse-equation}.
The rest of the assertion follows from Theorem~\ref{thm:ar-equation}.
\qed

\subsection{Morse inequalities.}
Recall that for an isolated invariant set $S$ we write $c_k(S):=\rank H^\kappa_k(S)$.
\begin{thm}
\label{thm:Morse-inequalities}
For a Morse decomposition $\cM$ of an isolated invariant set $S$ define
\[
   m_k(\cM):=\sum_{r\in{\scriptsize\PP}} c_k(M_r),
\]
where $c_k(M_r)=\rank H^\kappa_k(M_r)$.
Then  for any $k\in\ZZ^+$ we have the following inequalities
\begin{itemize}
   \item[(i)] (the strong Morse inequalities)
\[
   m_k(\cM)-m_{k-1}(\cM)+\cdots\pm m_0(\cM) \geq c_k(S)-c_{k-1}(S)+\cdots\pm c_0(S),
\]
   \item[(ii)] (the weak Morse inequalities)
\[
   m_k(\cM) \geq c_k(S).
\]
\end{itemize}
\end{thm}
\proof
Since $\cM$ is a Morse decomposition of $X$
also with respect to a linear extension of an admissible order, by a suitable
renaming of the elements of $\PP$ we may assume that $\PP=\{1,2,\ldots n\}$
with the order given by the order of integers.
Multiplying the equation \eqref{eq:Morse-equation} by the formal power series
\[
(1+t)^{-1}=1-t+t^2-t^3+\cdots
\]
and comparing the coefficients of both sides
we obtain the strong Morse Inequalities, because the polynomial $q$ is non-negative.
The weak Morse inequalities follow immediately from the strong Morse inequalities.
\qed


\section{Constructing combinatorial multivector fields from a cloud of vectors.}
\label{sec:algo}

In this section we present the algorithm used to generate some of the examples discussed in Section~\ref{sec:examples}.
It is a prototype algorithm intended to show the potential of the theory presented in this
paper in the algorithmic analysis of sampled dynamical systems
and in combinatorialization of flows given by differential equations.
The algorithm may be generalized to arbitrary polyhedral cellular complexes in $\RR^n$.
However, it is not necessarily optimal.
The question how to construct combinatorial multivector
fields from sampled flows in order to obtain a possibly best description of the dynamics
as well as how to compute the associated topological invariants efficiently is left for future investigation.

Consider $Z:=\setof{0,2,\ldots 2n}^2$. Let $E$ denote the set of open segments of length two with endpoints in $Z$
and let $S$ be the set of open squares of area four with vertices in $Z$. Then  $X:=Z\cup E\cup S$
is the collection of cells which makes the square  $[0,2n]^2$ a cellular complex.
We identify each cell with the coordinates of its center of mass. Thus, vertices have even coordinates,
squares have odd coordinates and edges have one coordinate even and one odd.
For an edge $e$ in $X$ we denote by $e^-,e^+$ its two vertices and we write $i^\bot(e)\in\{1,2\}$ for the index
of the axis orthogonal to the edge.

\begin{figure}

\begin{center}
  \includegraphics[width=0.45\textwidth]{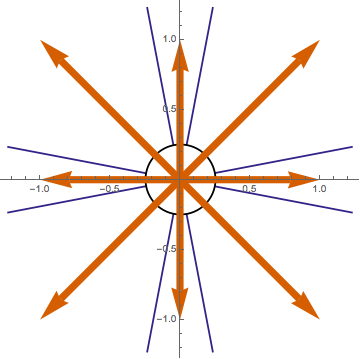}
\end{center}
  \caption{The normalization of the planar vectors via the division of the plane into nine regions: a disk around the origin, four cones
  along the four half-axes and four cones in the four quadrants of the plane.
  By a cone we mean here a cone excluding the disk at the origin.}
  \label{fig:nine-regions}
\end{figure}

Assume that we have a map $v:=(v_1,v_2):Z\to\RR^2$ which sends every vertex of $X$ to a planar, classical vector.
This may be a map which assigns to each point of $Z$ a vector in $\RR^2$
obtained from sampling a planar smooth vector field.
This may be also a map obtained from a physical experiment
or even selected randomly.

\begin{table}
\begin{tabbing}
boo\=boo\=boo\=boo\=boo\=boo\=boo\=boo\=\kill
\texttt{$\kwfunction$ CMVF($v$, $\mu$, $\epsilon$)}\\
\>\texttt{$\bar{v}:=D_{\mu,\epsilon}\circ v$;}\\
\>\texttt{\kwforeach $x\in X$ \kwdo $\theta(x):=x$ \kwenddo;}\\
\>\texttt{\kwforeach $x\in E$ \kwdo }\\
\>\>\texttt{$i:=i^\bot(x); s^-:=\bar{v}(x^-); s^+:=\bar{v}(x^+)$;}\\
\>\>\texttt{\kwif $s^-_i*s^+_i > 0$  \kwthen}\\
\>\>\>\texttt{$\theta(x):=$\kwif $s^-_1\neq 0$  \kwthen $x+(s^-_1,0)$ \kwelse $x+(0,s^-_2)$ \kwendif;}\\
\>\>\texttt{\kwendif;}\\
\>\texttt{\kwenddo;}\\
\>\texttt{\kwforeach $x\in Z$ \kwdo }\\
\>\>\texttt{$s:=\bar{v}(x)$; \kwif $s=0$ \kwcontinue;}\\
\>\>\texttt{\kwif $s_1*s_2=0$  \kwthen $\theta(x):=x+s$}\\
\>\>\texttt{\kwelse}\\
\>\>\>\texttt{ $x_1:=x+(s_1,0);\; x_2:=x+(0,s_2);$}\\
\>\>\>\texttt{ $t_1:=(\theta(x_1)=x_1);\; t_2:=(\theta(x_2)=x_2);$}\\
\>\>\>\texttt{ $c_1:=t_1\kwand \bar{v}_1(x_1+s_1)*s_1<0;$}\\
\>\>\>\texttt{ $c_2:=t_2\kwand \bar{v}_2(x_2+s_2)*s_2<0;$}\\
\>\>\>\texttt{\kwif $c_1$ \kwand $c_2$  \kwthen \kwcontinue;}\\
\>\>\>\texttt{\kwif \kwnot $c_1$ \kwand \kwnot $c_2$  \kwthen }\\
\>\>\>\>\texttt{ $\theta(x_1):=\theta(x_2):=\theta(x):=x+s$ ;}\\
\>\>\>\texttt{\kwelse }\\
\>\>\>\>\texttt{\kwif $c_1$ \kwand $t_2$ \kwthen $\theta(x):=x_1$ \kwendif;}\\
\>\>\>\>\texttt{\kwif $c_2$ \kwand $t_1$ \kwthen $\theta(x):=x_2$ \kwendif;}\\
\>\>\>\texttt{\kwendif;}\\
\>\>\texttt{\kwendif;}\\
\>\texttt{\kwenddo;}\\
\>\texttt{\kwreturn $\theta$;}\\
\texttt{\kw{endfunction};}\\
\end{tabbing}
\caption{Algorithm CMVF constructing a combinatorial multivector field from a collection of vectors
on an integer lattice in the plane. }
\label{tab:cmv-alg}
\end{table}

Let $\sgn:\RR\to\{-1,0,1\}$ denote the sign function.
Consider the map $\Arg:\RR^2\setminus\{0\}\to [-\pi,\pi)$ given by
\[
   \Arg(x_1,x_2):=\begin{cases}
                    \arccos \frac{x_1}{\sqrt{x_1^2+x_2^2}} & \text{if $x_2\geq 0$, }\\
                    -\arccos \frac{x_1}{\sqrt{x_1^2+x_2^2}} & \text{if $x_2< 0$. }
                  \end{cases}
\]
and the planar map $D_{\mu,\epsilon}:\RR^2\to\RR^2$ given by
\[
    D_{\mu,\epsilon}(x_1,x_2):=
    \begin{cases}
      0            & \text{if $x_1^2+x_2^2\leq \epsilon^2$,}\\
      (\sgn x_1,0) & \text{if $|\Arg x|\leq\mu$ or $|\Arg x|\geq\pi-\mu$,}\\
      (0,\sgn x_2) & \text{if $|\Arg x-\pi/2|\leq\mu$}\\
                   & \text{or $|\Arg x+\pi/2|\leq\mu$,}\\
      (\sgn x_1,\sgn x_2) & \text{otherwise.}
    \end{cases}
\]
The map $D_{\mu,\epsilon}$ normalizes the vectors in the plane in the sense that it sends
each planar vector $v$ to one of the nine vectors with coordinates in $\{-1,0,1\}$
depending on the location of $v$ in one of the nine regions marked in Figure~\ref{fig:nine-regions}.

The algorithm is presented in Table~\ref{tab:cmv-alg}.
Its rough description is as follows. It accepts on input a triple $(v,\mu,\epsilon)$,
where $v:Z\to\RR^2$ is a map, $\mu\in [0,\pi/4]$ and $\epsilon>0$.
The map $v$ is normalized to $\bar{v}$ by applying  the map  $D_{\mu,\epsilon}$.
The algorithm defines a map $\theta:X\to X$ such that $\cV_\theta:=\setof{\theta^{-1}(y)\mid y\in\im\theta}$
encodes a combinatorial multivector field via Theorem~\ref{thm:theta}.
The construction proceeds in three steps. In the first step the map is initialized to identity.
This corresponds to a combinatorial multivector field consisting only of singletons.
In the second step each edge is tested whether the normalized vectors at its endpoints
projected to the axis perpendicular to the edge coincide.
If this is the case, the edge is paired with the square pointed by both projections.
In the third step each vertex $z\in Z$ is analysed. If the normalized vector in $\bar{v}(z)$
is parallel to one of the axes, the vertex $z$ is paired with the corresponding edge.
Otherwise, it is combined into one multivector with the corresponding square and adjacent edges
unless this would generate a conflict with one of the neighbouring vertices.

Note that the only case when strict multivectors may be generated is in the third step when
the normalized vector is not parallel to any axis. The chances that this happens decrease
when the parameter $\mu$ is increased. The normalized vector is always parallel to one of the axes
when $\mu\geq \pi/4$. Thus, if $\mu\geq \pi/4$, then the algorithm returns a combinatorial vector
field, that is a combinatorial multivector field with no strict multivectors.

\begin{thm}
\label{thm:algo}
Consider the algorithm in Table~\ref{tab:cmv-alg}.
When applied to a map $v$ and  numbers $\mu> 0$, $\epsilon>0$
it always stops returning a map $\theta:X\to X$ which satisfies conditions (i)-(iii)
of Theorem~\ref{thm:theta}. Thus, the collection $\setof{\theta^{-1}(y)\mid y\in\im \theta}$
is a combinatorial multivector field on $X$.
\end{thm}
\proof
  After completing the first $\kwforeach$ loop the map $\theta$, as the identity,
trivially satisfies conditions (i)-(iii) of Theorem~\ref{thm:theta}.
Thus, it suffices to observe that each modification of $\theta$ in the subsequent loops
preserves theses properties. One can easily check that this is indeed the case.
\qed

\section{Extensions.}
\label{sec:extensions}

In this section we briefly indicate the possibilities
of some extensions of the theory presented in this paper.
The details will be presented elsewhere.

\begin{figure}
\begin{center}
  \hspace*{18pt}\includegraphics[width=0.65\textwidth]{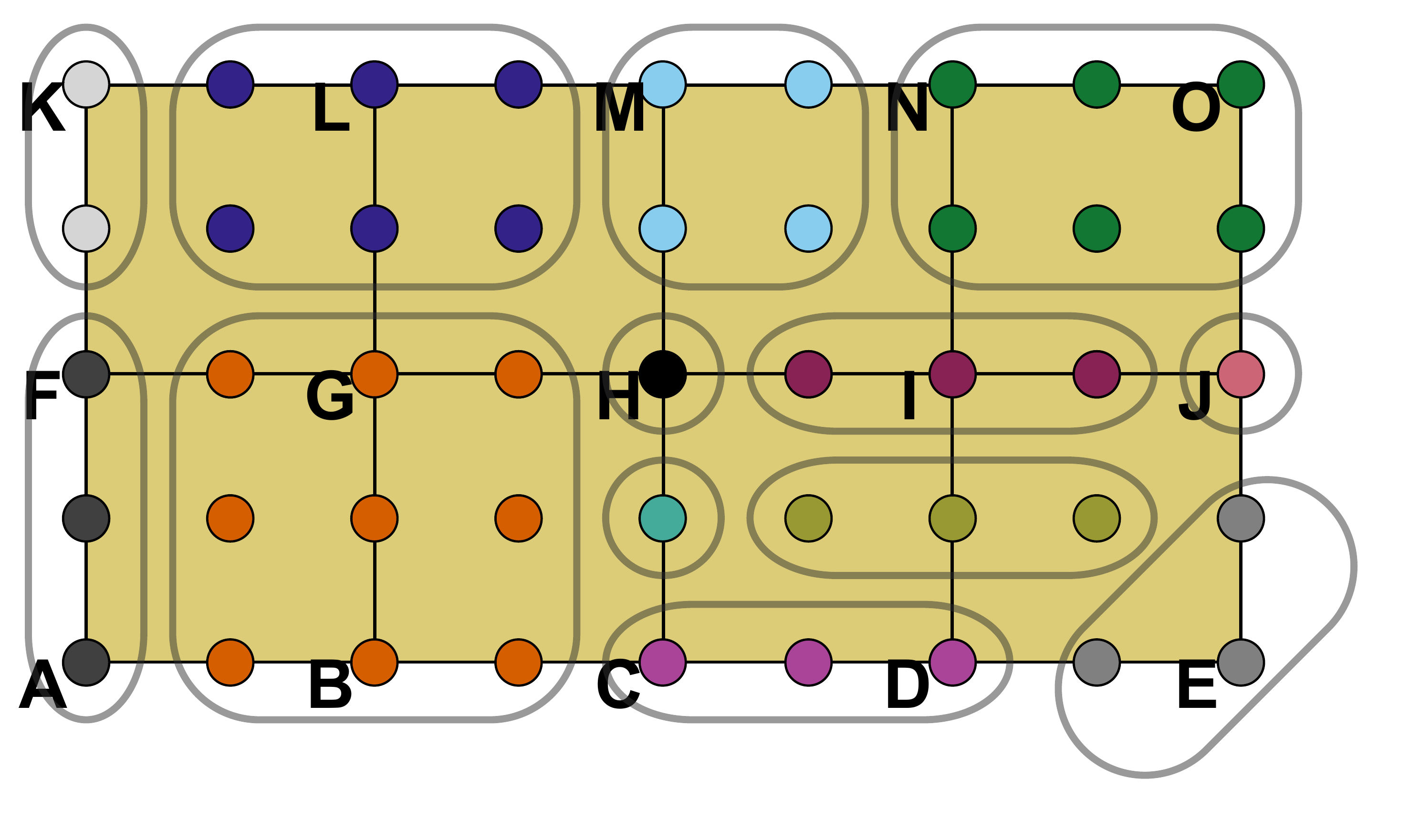}\\
  \includegraphics[width=0.6\textwidth]{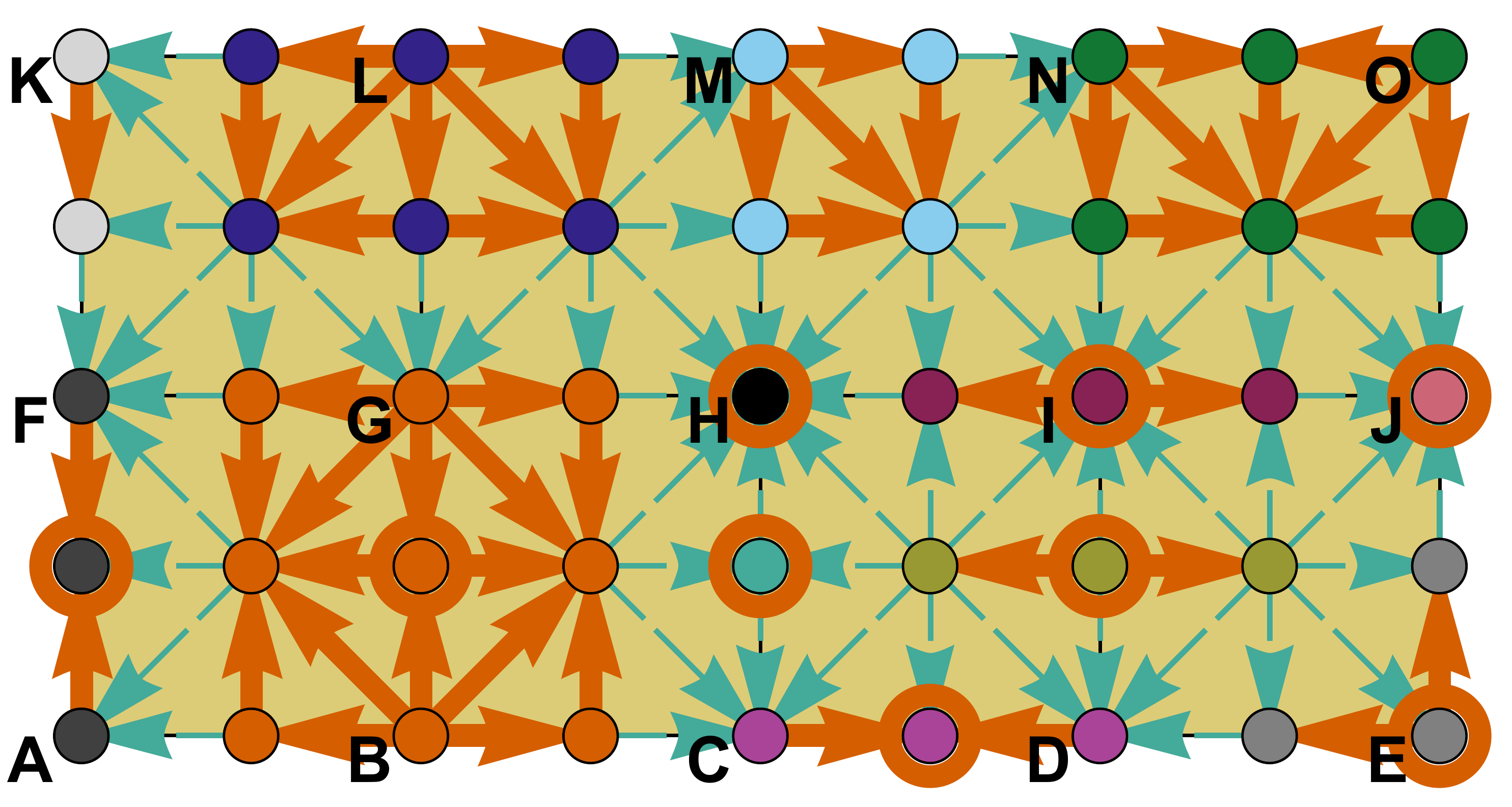}
\end{center}
  \caption{A generalized multivector field as a partition of a Lefschetz complex (top)
  and the associated  directed graph $G_\cV$ (bottom).
  The up-arrows and loops are marked by thick solid lines.
  The down-arrows are marked by thin dashed lines.
  Only one loop for each critical vector is marked.
  Multivectors in the upper row of squares are regular.
  Multivectors in the lower row of squares are critical.
  }
  \label{fig:genMVFpartition}
\end{figure}

\subsection{Generalized multivector fields.}
The requirement that a multivector has a unique maximal element with respect
to the partial order $\leq_\kappa$ is convenient. It simplifies the analysis of
combinatorial multivector fields, because every non-dominant cell $x$ admits
precisely one arrow in $G_\cV$ originating in $x$, namely the up-arrow from $x$ to
$x^\star$. It is also not very restrictive in combinatorial modelling of a dynamical
system. But, in some situations the requirement may be inconvenient. For instance, if sampling is performed
near a hyperbolic or repelling fixed point, it may happen that the sampled vectors point into
different top dimensional cells. Or, when a parameterized family of dynamical systems is studied
(see Section~\ref{sec:continuation}),
it may be natural to model a collection of nearby dynamical systems
by one combinatorial multivector field. In such a situation the uniqueness
requirement will not be satisfied in places where neighbouring systems point into different cells.

Taking these limitations into account it is reasonable to consider
the following extension of the concept of a combinatorial multivector field.
A {\em generalized multivector} in a Lefschetz complex $(X,\kappa)$ is a proper collection of cells in $X$.
A generalized multivector $V$ is {\em critical} if $H^\kappa(V)\neq 0$.
A {\em generalized multivector field} $\cV$ is a partition of $X$ into generalized multivectors.
For a cell $x\in X$ we denote by $[x]$ the unique generalized multivector in $\cV$ to which $x$ belongs.
The associated generalized directed graph $G_\cV$ has vertices in $X$ and an arrow from $x$ to $y$ if one of the following
conditions is satisfied
\begin{eqnarray}
  && \dim x < \dim y  \text{ and } y\in [x] \text{ (an {\em up-arrow}),}\label{eq:gadhr1}\\
  && \dim x > \dim y  \text{ and } y\in\mouth[x]   \text{ (a {\em down-arrow}),}\label{eq:gadhr2}\\
  && x=y \text{ and } [y] \text{ is critical (a {\em loop}).}\label{eq:gadhr3}
\end{eqnarray}
We write $y \adhl_{\cV} x$ if there is an arrow from $x$ to $y$ in $G_\cV$,
denote by $\leqV$ the preorder induced by $\adhl_{\cV}$, interpret $\adhl_{\cV}$
as a multivalued map $\Pi_{\cV}: X\mvmap X$ and study the associated dynamics.
An example of a generalized multivector field together with the associated up-arrows
is presented in Figure~\ref{fig:genMVFpartition}.

The main ideas of the theory presented in this paper extend to this generalized case.
The proofs get more complicated, because there are more cases to analyse. This is caused
by the fact that more than one up-arrow may originate from a given cell.

It is tempting
to simplify the analysis by gluing vertices of $G_\cV$ belonging to the same multivector.
But, when proceeding this way, the original phase space is lost.
And, the concept of index pair makes sense only in the original space, because in general
the sets $P_1$, $P_2$ in an index pair $P=(P_1,P_2)$ are not $\cV$-compatible.
Also, since in the gluing process two different combinatorial multivector fields would modify the space in a different way,
it would not be possible to compare their dynamics.
In consequence, the  concepts of perturbation and continuation briefly explained in the following section
would not make sense.

\begin{figure}
\begin{center}
  \includegraphics[width=0.57\textwidth]{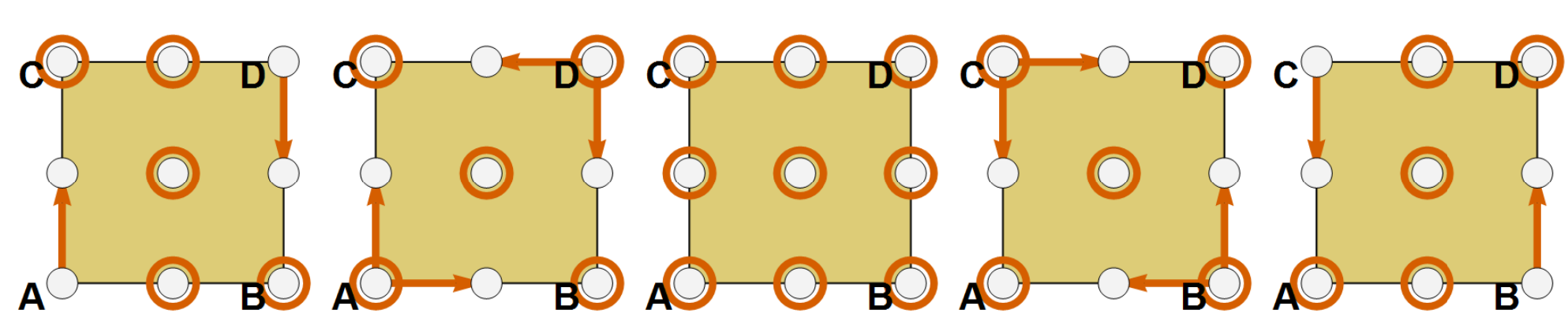}\\
  \includegraphics[width=0.8\textwidth]{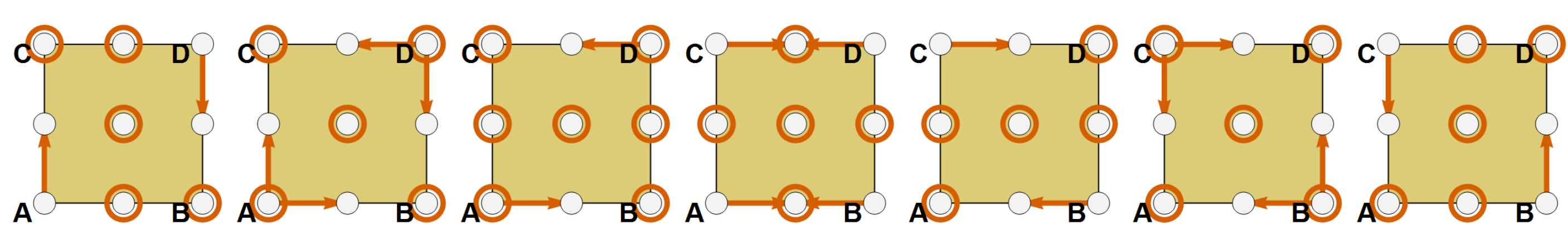}
\end{center}
  \caption{Two parameterized families of combinatorial multivector fields: $\cV_1$, $\cV_2$, $\cV_3$, $\cV_4$, $\cV_5$ (top from left to right)
  and $\cW_1$, $\cW_2$, $\cW_3$, $\cW_4$, $\cW_5$, $\cW_6$, $\cW_7$ (bottom from left to right). Note that
  }
  \label{fig:continuation}
\end{figure}

\subsection{Perturbations and continuation.}
\label{sec:continuation}
Among the fundamental features of the classical Conley index theory is the continuation property \cite[Chapter IV]{Co78}.
An analogous property may be formulated in the combinatorial case.
Let $\cV$, $\bar{\cV}$ be two combinatorial multivector fields on a Lefschetz complex $X$.
We say that $\bar{\cV}$ is a {\em perturbation} of $\cV$ if $\bar{\cV}$ is a refinement of $\cV$ or $\cV$
is a refinement of $\bar{\cV}$.
A {\em parameterized family} of combinatorial multivector fields is
a sequence $\cV_1,\cV_2,\ldots,\cV_n$ of combinatorial multivector fields on $X$
such that $\cV_{i+1}$ is a perturbation of $\cV_i$.
Assume $S$ is an isolated invariant set of $\cV$ and $\bar{S}$ is an isolated invariant set of $\bar{\cV}$.
We say that $S$ and $\bar{S}$ are {\em related by a direct continuation} if there is a pair $P=(P_1,P_2)$ of closed sets
in $X$ such that $P$ is an index pair for $S$ with respect to $\cV$ and for $\bar{S}$ with respect to $\bar{\cV}$.
We say that $S$ and $\bar{S}$ are {\em related by continuation} along
a parameterized family $\cV=\cV_1,\cV_2,\ldots,\cV_n=\bar{\cV}$
if  there exists isolated invariant sets $S_i$ of $\cV_i$
such that $S_1=S$, $S_n=\bar{S}$ and  $S_i$, $S_{i+1}$ are
related by a direct continuation.
 In a similar way one can define continuation of Morse decompositions of isolated invariant sets.

Figure~\ref{fig:continuation} presents two parameterized families of combinatorial multivector fields:
$\cV_1$, $\cV_2$, $\cV_3$, $\cV_4$, $\cV_5$ in the top row and
$\cW_1$, $\cW_2$, $\cW_3$, $\cW_4$, $\cW_5$, $\cW_6$, $\cW_7$
in the bottom row. Note that $\cV_1=\cW_1$ and $\cV_5=\cW_7$.
The singleton $\{AB\}$ is an isolated invariant set of $\cV_1$.
The set $\{A,AB,AC\}$ is an isolated invariant set of $\cV_2$.
These two sets are related by a direct continuation, because
the pair $(\{A,B,C,AB,AC\},\{B,C\})$ is an index pair
for $\{AB\}$ with respect to $\cV_1$ and for $\{A,AB,AC\}$ with respect to $\cV_2$.
Note that $\{AB\}$ is also an isolated invariant set of $\cV_5$ but it is not related
by continuation along $\cV_1$, $\cV_2$, $\cV_3$, $\cV_4$, $\cV_5$ to  $\{AB\}$ as an isolated invariant set of $\cV_1$.
But $\{AB\}$ as an isolated invariant set for $\cW_1$ is related by continuation
along $\cW_1$, $\cW_2$, $\cW_3$, $\cW_4$, $\cW_5$, $\cW_6$, $\cW_7$
to $\{CD\}$ as an isolated invariant set for $\cW_7$. Actually, the Morse decomposition
of $\cW_1$ consisting of Morse sets $\{ABCD\}$, $\{AB\}$, $\{CD\}$, $\{B\}$, $\{C\}$ is related by continuation
along $\cW_1$, $\cW_2$, $\cW_3$, $\cW_4$, $\cW_5$, $\cW_6$, $\cW_7$ to the Morse decomposition
of $\cW_7$ consisting of Morse sets $\{ABCD\}$, $\{CD\}$, $\{AB\}$, $\{C\}$, $\{B\}$.
In particular, the isolated invariant sets $\{ABCD\}$, $\{AB\}$, $\{CD\}$, $\{B\}$, $\{C\}$ of $\cW_1$
are related by continuation respectively to the isolated invariant sets
$\{ABCD\}$, $\{CD\}$, $\{AB\}$, $\{C\}$, $\{B\}$ of $\cW_7$.

\subsection{Homotopy Conley index for combinatorial multivector fields on CW complexes.}
Consider a finite regular CW complex $X$ with the collection of cells $K$.
We say that $V\subset K$ is a multivector if $|V|$, the union of all cells in $V$, is a proper subset of $X$.
We say that a multivector $V$ is {\em regular} if $\mouth |V|:=\cl|V|\setminus|V|$ is a deformation retract of $\cl |V|\subset X$.
Otherwise we call $V$ {\em critical}. We define the {\em generalized multivector field} $\cV$ on $X$
as a partition of $K$ into generalized multivectors. Using the modified concepts of regular and critical
multivectors as in the case of Lefschetz complexes we introduce the directed graph $G_\cV$ and the associated dynamics.
This lets us introduce the concept of an isolated invariant set $S$ and the associated Conley index defined as the homotopy
type of $P_1/P_2$ for an index pair $(P_1,P_2)$ isolating $S$. The proofs rely on gluing the homotopies provided by the
deformation retractions of the individual regular multivectors.

\subsection{Discrete Morse theory for combinatorial multivector fields.}
The theory presented in this paper is built on the concept of a combinatorial
multivector field, an extension of Forman's combinatorial vector field used to construct
the discrete Morse theory. A natural question is whether one can use multivectors instead of vectors
in this construction under the additional assumption that  all the critical cells are non-degenerate.
It seems that the collapsing approach to discrete Morse theory presented in \cite{BrElJuMr2015}
extends to the case of combinatorial multivector
fields on CW complexes defined in the preceding section. In the algebraic case of Lefschetz
complexes one needs an extra assumption that the closure of every regular multivector
is chain homotopy equivalent to its mouth.

\section{Conclusions and future research.}
\label{sec:conclu}

The presented theory shows that combinatorialization of dynamics, started by Forman's paper
\cite{Fo98b}, may be extended to cover such concepts as isolated invariant set, Conley
index, attractor, repeller, attractor-repeller pair and Morse decomposition.
Moreover, the broader concept
of a combinatorial multivector field introduced in the present paper seems to be more flexible
than the original Forman's concept of combinatorial vector field and shall serve better
the needs of the combinatorialization of classical topological dynamics.
In particular, some classical concepts in dynamics, for instance the homoclinic connection,
do not have counterparts in the theory of combinatorial vector fields but do have
in the case of combinatorial multivector fields.

The theory proposed in this paper is far from being completed.
Besides the extensions discussed in the previous section there are several directions of research to be undertaken.
For instance,  connection matrices constitute an essential part of the Conley index theory.
It is natural to expect that they have some counterpart in the combinatorial setting.
Several authors proposed a generalization of the Conley index theory to the case of time-discrete dynamical systems
(iterates of maps) \cite{robbin-salamon-1985-1988,mrozek-1990a,Szymczak1995,FranksRicheson2000}.
A challenging question is what is a combinatorial analogue of a time-discrete
dynamical system and how to construct Conley index theory for such systems.
Poincar\'e maps for combinatorial  multivector fields
may be a natural place to start investigations in this direction.
Also, it would be interesting to understand the relation between the dynamics of a combinatorial multivector field and its Forman refinements.

The real significance of combinatorial multivector fields will be seen
only after constructing bridges between combinatorial and classical dynamics.
Some work towards a bridge from combinatorial to classical dynamics has been undertaken in \cite{KMW}.
However, from the point of view of applications a reverse approach is more important:
constructing a combinatorial multivector field modelling a classical flow.
The algorithm presented in this paper is only a hint that such constructions are possible.
Obviously, finite resolution of the combinatorial setting prohibits any one-to-one
correspondence, but approximation schemes representing the classical
flow up to the resolution controlled by the size of the cells in the approximation should be possible.
Depending on the goal, there are at least two options. If the goal is the rigourous
description of the classical dynamics by means of a combinatorial multivector field,
the techniques developed in \cite{BoKaMi2007} might be of interest.
If rigour is not necessary, the use of a piecewise constant approximation
of a vector field \cite{Szy2012} to construct a combinatorial multivector field
is an interesting option.

\begin{center}
{\bf Acknowledgments.}
\end{center}
I am indebted to Abhishek Rathod who presented some ideas of Conley theory
in the setting of Forman's combinatorial vector fields in a talk presented at the DyToComp 2012
conference in Bedlewo, Poland. This was a direct inspiration to undertake the present research.
I am grateful to Ulrich Bauer and Herbert Edelsbrunner for pointing out to me
the concept of the generalized Morse matchings.
This conversation lead directly to the idea of the combinatorial multivector field
presented in this paper. Also, I would like to thank Tamal Dey
for several valuable comments which let me enhance the presentation of the paper.

\end{document}